\input amstex
\documentstyle{amsppt}
\document
\input xypic
\NoBlackBoxes

\define\A{\Bbb A}
\define\BB{\Bbb B}
\define\J{\Bbb J}
\define\Q{\Bbb Q}
\define\G{\Bbb G}
\define\Z{\Bbb Z}

\define\F{\Bbb F}
\define\R{\Bbb R}
\define\C{\Bbb C}

\define\SS{\Bbb S}
\define\T{\Bbb T}
\define\X{\Bbb X}
\define\Y{\Bbb Y}

\define\cal{\Cal}

\topmatter
\rightheadtext{The image of {\sl l}-adic Galois representations,}
\title On the image of {\sl l}-adic Galois
representations for abelian
varieties of type I and II
\endtitle
\author G. Banaszak, W. Gajda, P. Kraso\'n
\endauthor
\address
Department of Mathematics, Adam Mickiewicz University,
Pozna\'{n}, Poland
\endaddress
\email banaszak\@math.amu.edu.pl
\endemail
\address
Department of Mathematics, Adam Mickiewicz University,
Pozna\'{n}, Poland
\endaddress
\email gajda\@math.amu.edu.pl
\endemail
\address
current: Centre de Recerca Matematica, Bellaterra, Spain
\endaddress
\address
Department of Mathematics, Szczecin University,
Szczecin, Poland
\endaddress
\email
krason\@sus.univ.szczecin.pl
\endemail

\abstract In this paper we investigate the image of the $l$-adic
representation attached to the Tate module of an abelian variety
over a number field with endomorphism algebra of type I or II in
the Albert classification. We compute the image explicitly  and
verify the classical conjectures of Mumford-Tate, Hodge, Lang and
Tate, for a large family of abelian varieties of type I and II. In
addition, for this family, we prove an analogue of the open image
theorem of Serre.
\endabstract
\endtopmatter

\subhead  1. Introduction
\endsubhead

Let $A$ be an abelian variety defined over a number field $F.$ Let
$l$ be an odd prime number. In this paper we study the images of
the $l$-adic representation $\rho_{l}: G_F \longrightarrow
GL(T_{l}(A))$ and the $mod\,\, l$ representation ${\overline
\rho_{l}}: G_F \longrightarrow GL(A[l])$ of the absolute Galois
group $G_F=G({\bar F}/F)$ of the field $F,$ associated with the
Tate module, for $A$ of type I or II in the Albert classification
list cf. [M]. In our previous paper on the subject cf. [BGK], we
computed the images of the Galois representations for some abelian
varieties with real (type I) and complex multiplications (type IV)
by the field $E{=}End_F(A)\otimes \Q$ and for $l$ which splits
completely in the field $E$ {\sl loc. cit.}, Theorem 2.1 and
Theorem 5.3.
\medskip

In the present paper we extend results proven in [BGK] to a larger
class (cf. Definition of class ${\cal A}$ below) of abelian
varieties which includes varieties with non-commutative algebras of
endomorphisms, and to almost all prime numbers $l.$ In order to get
these results, we had to implement the Weil restriction functor
$R_{L/K}$ for a finite extension of fields $L/K.$ In section 2 of
the paper we give an explicit description of the Weil restriction
functor for affine group schemes which we use in the following
sections. In a very short section 3 we prove two general lemmas
about bilinear forms which we apply to Weil pairing in the following
section. Further in section 4, we collect some auxiliary facts about
abelian varieties. In section 5 we obtain the integral versions of
the results of Chi cf. [C2], for abelian varieties of type II and
compute Lie algebras and endomorphism algebras corresponding to the
$\lambda$-adic
 representations related to the Tate module of $A.$
In section 6 we prove the main results of the paper concerning the
image of Galois representation $\rho_l$ and the Zariski closures
$G_{l}^{alg},$ $G(l)^{alg}$ of the images of the
$\Q_l$-representation $\rho_{l}\otimes \Q_l : G_F \rightarrow
GL(V_{l}(A))$ and of the mod $l$-representation ${\overline
\rho_{l}} : G_F \rightarrow GL(A[l]).$
\medskip

Main results proven in this paper concern the following class of
abelian varieties:

\proclaim{Definition of class ${\cal A}$}

\noindent We say that an abelian variety $A/F,$ defined over a
number field $F$ is of class ${\cal A},$ if the following conditions
hold:

\roster \item"{(i)}" $A$ is a simple, principally polarized abelian
variety of dimension $g$
 \item"{(ii)}" ${\cal R}=End_{\bar F} (A)=End_{F}(A)$ and the endomorphism algebra
 $D = {\cal R} \otimes_{\Z} \Q,$ is of type I or II in the Albert
list of  division algebras with involution  (cf. [Mu], p. 210).
 \item"{(iii)}" the field $F$ is such that for
every $l$ the Zariski closure $G_l^{alg}$ of $\rho_l(G_F)$ in
$GL_{2g}/\Q_l$ is a connected algebraic group
 \item"{(iv)}" $ g = h e d,$ where $h$ is an odd integer,
 $e = [E: Q]$ is the degree of the center $E$ of $D$ and
 $d^2 = [D:\, E].$
\endroster
\endproclaim

\noindent We have chosen to work with principal polarizations,
however main results of this paper have their analogs for any simple
abelian variety $A$ with a fixed polarization, provided $A$
satisfies the above conditions (ii), (iii) and (iv).
 The most restrictive of the conditions in the definition of class $\Cal A$
 is condition (iv) on the dimension of the variety $A.$
 We need this condition to perform  computations with Lie algebras in the proof
 of Lemma 5.33, which are based on
 an application of the miniscule conjecture cf. [P].
Note that due to results of Serre,
 the assumption (iii) is not very restrictive. It follows by [Se1] and
[Se4], that for an abelian variety $A$ defined over a number field
$K,$ there exists a finite extension $K^{conn}/K$ for which the
Zariski closure of the group $\rho_l(G_{K^{conn}})$ in $GL$ is a
connected variety for any prime $l$. Hence, to make $A$ meet the
 condition (iii), it is enough to enlarge the base field, if necessary.
 Note that the field $K^{conn}$ can
 be determined in purely algebraic terms, as the intersection of a family of 
fields
 of division points on the abelian variety $A$ cf. [LP2], Therem 0.1.
 \bigskip\bigskip

\noindent {\bf Main results}
\medskip

\proclaim{Theorem A} [Theorem 6.9]

\noindent If $A$ is an abelian variety of class ${\cal A},$ then for
$l \gg 0,$ we have equalities of group schemes:
$$(G_{l}^{alg})^{\prime} \,\, = \,\, \prod_{\lambda | l} \,\,R_{E_{\lambda}/\Q_l}
(Sp_{2h})$$
$$(G(l)^{alg})^{\prime} \,\, = \,\, \prod_{\lambda | l} \,\,R_{k_{\lambda}/\F_l}
(Sp_{2h}),$$ where $G^{\prime}$ stands for the commutator subgroup
of an algebraic group $G,$ and $R_{L/K}(-)$ denotes the Weil
restriction functor.
\endproclaim
\medskip

\proclaim{Theorem B} [Theorem 6.16]

\noindent If $A$ is an abelian variety of class ${\cal A},$  then
for $l \gg 0,$ we have:

$${\overline{\rho_l}}(G_F^{\prime})
\,\, = \,\, \prod_{\lambda | l} \,\, Sp_{2h}(k_{\lambda})\,\, =\,\,
Sp_{2h}({\cal O}_{E}/l {\cal O}_{E} )$$
$$\rho_l\bigl(\,\overline{G_F^{\prime}}\,\bigr)
\,\, = \,\, \prod_{\lambda | l} \,\, Sp_{2h}({\cal O}_{\lambda})
\,\, =\,\, Sp_{2h}({\cal O}_{E} \otimes_{\Z} \Z_l),$$ where
$\overline{G_F^{\prime}}$ is the closure of $G_F^{\prime}$ in the
profinite topology in $G_F.$
\endproclaim

As an application of Theorem A we obtain: \proclaim{Theorem C}
[Theorem 7.12]

\noindent \noindent If $A$ is an abelian variety of class ${\cal
A},$ then
$$G_l^{alg}= MT(A)\otimes \Q_l,$$
for every prime number $l$, where $MT(A)$ denotes the Mumford-Tate
group of $A,$ i.e., the Mumford -Tate conjecture holds true for $A.$
\endproclaim
\medskip

Using the approach initiated by Tankeev [Ta5] and Ribet [R2], futher
developed by V.K. Murty  [Mu] combined with some extra work on the
Hodge groups in section 7, we obtain:
\medskip

\proclaim{Theorem D} [Thorems 7.34, 7.35]

\noindent If $A$ is an abelian variety of class ${\cal A},$ then the
Hodge conjecture and the Tate conjecture on the algebraic cycle
maps, hold true for the abelian variety $A.$
\endproclaim
\medskip

\noindent In the past there has been an extensive work on
the Mumford-Tate, Tate and Hodge conjectures for abelian varieties.
Some special cases of the conjectures were verified for other
classes of abelian varieties, see for example papers: [Ab], [C1], [C2],
[Mu], [P], [Po], [R2], [Se1], [Se5], [Ta1], [Ta2], [Ta3], [Ta4] and
[BGK]. For more complete list of results concerning the Hodge conjecture
see [G].
\medskip

Moreover, using a result of Wintenberger (cf. [Wi], Cor. 1, p.5), we
were able to verify that for $A$ of class ${\cal A},$ the group
$\rho_{l}(G_F)$ contains the group of all the homotheties in
$GL_{T_l(A)}(\Z_l)$ for  $l \gg 0,$ i.e., the Lang conjecture
holds true for $A$ cf. Theorem 7.38.
\medskip

As a final application of the method developed in this paper, we
prove an analogue of the open image theorem of Serre cf. [Se1] for
the class of abelian varieties with which we work.
\medskip

\proclaim{Theorem E} [Theorem 7.42]

\noindent If $A$ is an abelian variety of class ${\cal A},$ then
for every prime number $l,$ the image $\rho_{l}(G_F)$ is open in
the group $C_{{\cal R}}(GSp_{(\Lambda,\, \psi)})(\Z_l)$ of
$\Z_l$-points of the commutant of ${\cal R}{=}End\,A$ in the group
$GSp_{(\Lambda,\, \psi)}$ of symplectic similitudes of the
bilinear form $\psi:\Lambda\times \Lambda\longrightarrow \Z$
associated with the polarization of $A.$ In addition, for $l \gg
0$ we have:
$$\rho_{l}(\overline{G_F^{\prime}}) \,\, =\,\,
C_{{\cal R}}(Sp_{(\Lambda,\, \psi)})(\Z_l).$$
\endproclaim
\medskip

\noindent As an immediate corollary of Theorem E we obtain that
for any $A$ of class ${\cal A}$ and for every $l,$ the group
$\rho_{l}(G_F)$ is open in ${\cal G}_{l}^{alg}(\Z_l)$ (in the
$l$-adic topology), where ${\cal G}_{l}^{alg}$ is the Zariski
closure of $\rho_l(G_F)$ in $GL_{2g}/\Z_l.$ cf. Theorem 7.48.
\bigskip\bigskip

\noindent \subhead  2. Weil restriction functor
$R_{E/K}$ for affine schemes and Lie algebras
\endsubhead

In this section we describe Weil restriction functor and its basic
properties which will be used in the paper c.f. [V1], [V2,\, pp. 37-40],
[W1] and [W2,\, pp. 4-9]. Let $E/K$
be a separable field extension of degree $n.$ Let
$\{\sigma_1, \sigma_2,\dots , \sigma_n\}$ denote the set of all
imbeddings $E \rightarrow E^{\sigma_i} \subset {\overline K}$
fixing $K.$ Define $M$ to be the composite of the fields
$E^{\sigma_i}$ $$M = E^{\sigma_1}\dots E^{\sigma_n}.$$
Let $X =[x_1, x_2,\dots x_r]$ denote a multivariable. For polynomials $f_k = f_k(X) \in E[X],$
$1\leq k\leq s,$ we denote by $I = (f_1, f_2, \dots, f_s)$ the ideal generated by the
$f_k $'s  and put $I^{\sigma_i} = (f_{1}^{\sigma_i}(X), f_{2}^{\sigma_i}(X), \dots,
f_{s}^{\sigma_i}(X))$ for any $1\leq i\leq n.$  Let
$A = E[X]/ I.$  Define $E$-algebras $A^{\sigma_i}$ and ${\overline A}$ as follows:
$$A^{\sigma_i} = A \otimes_{E, \sigma_i} M\, \cong\, M[X]/\, I^{\sigma_{i}}M[X],$$
$${\overline A} = A^{\sigma_1}\otimes_{M}\dots\otimes_{M}  A^{\sigma_n}.$$
Let $X^{\sigma_1},\dots, X^{\sigma_n}$ denote the multivariables
$$X^{\sigma_i} = [x_{i,1}, x_{i,2},\dots, x_{i,r}]$$ on which the Galois group
$G = G(M/K)$ acts naturally on the right in the following way: since $E = K(\gamma),$
for any imbedding $\sigma_i$ and any $\sigma \in G$ we get
$(\gamma^{\sigma_{i}})^{\sigma} = \gamma^{\sigma_{j}}$ for some $1\leq j\leq n.$ We define:
$$(X^{\sigma_{i}})^{\sigma} = X^{\sigma_{j}}.$$
We see that
$${\overline A} \,\cong\, M[X^{\sigma_1},\dots, X^{\sigma_n}]/\, (I_1 + \dots + I_n),$$
where $I_k = M[X^{\sigma_1},\dots, X^{\sigma_n}]I_{(k)}$
and $I_{(k)} = (f_{1}^{\sigma_k}(X^{\sigma_k}),\dots, f_{s}^{\sigma_k}(X^{\sigma_k})),$
for any $1\leq k\leq n.$
\medskip

\proclaim{Lemma 2.1}
$${{\overline A}}^G \otimes_{K} M \cong {\overline A}.$$
\endproclaim

\demo{\sl Proof}
Let $\alpha_1,\dots,\alpha_n$ be a basis of $E$ over $K.$ It is clear that
$$\sum_{i=1}^n \alpha_{j}^{\sigma_{i}} X^{\sigma_{i}} \,\,\in\,\, {{\overline A}}^G.$$
Since $[\alpha_{j}^{\sigma_{i}}]_{i,j}$ is an invertible matrix with coefficients in $M,$
we observe that $X^{\sigma_{1}}, \dots, X^{\sigma_{n}}$ are in the subalgebra of
${\overline A}$ generated by $M$ and ${{\overline A}}^G.$ But
$X^{\sigma_{1}}, \dots, X^{\sigma_{n}}$ and $M$ generate ${\overline A}$ as an algebra.
\qed\enddemo

\remark{\bf Remark 2.2} Notice that the elements $\sum_{i=1}^n
\alpha_{j}^{\sigma_{i}} X^{\sigma_{i}}$ for $j = 1, \dots, n$
generate ${{\overline A}}^G$ as a K-algebra. Indeed if $C$
denotes the $K$-subalgebra of ${{\overline A}}^G$ generated by
these elements and if $C$ were smaller than ${{\overline A}}^G,$
then $C \otimes_K M$ would be smaller then  ${{\overline A}}^G
\otimes_K M,$ contrary to Lemma 2.1.
\endremark

\medskip
\proclaim{Definition 2.3} Put $V = spec\, A,$ and $W = spec\, {{\overline A}}^G.$
Weil's restriction functor $R_{E/K}$ is defined by the following formula:
$$R_{E/K}(V) = W.$$
\endproclaim
\medskip

\noindent
Note that we have the following isomorphisms:
$$W\otimes_{K} M \,\,
=\,\, spec\, ({{\overline A}}^G \otimes_{K} M) \,\, \cong\,\,
spec\, {\overline A} \cong$$
$$ spec\, (A^{\sigma_1}\otimes_{M}\dots\otimes_{M}  A^{\sigma_n}) \,\, \cong \,\, (V\otimes_{E, \sigma_{1}} M) \otimes_{M}\dots \otimes_{M} (V\otimes_{E, \sigma_{n}} M),$$
hence
$$R_{E/K}(V)\otimes_{K} M \cong (V\otimes_{E, \sigma_{1}} M) \otimes_{M}\dots \otimes_{M}
(V\otimes_{E, \sigma_{n}} M).$$

\medskip
\proclaim{Lemma 2.4} Let $V^{\prime} \subset V$ be a closed
imbedding of affine schemes over $E.$ Then
$R_{E/K}(V^{\prime}) \subset R_{E/K}(V)$ is a closed
imbedding of affine schemes over~$K.$
\endproclaim
\demo{\sl Proof} We can assume that $V = spec\, (E[X]/I)$ and
$V^{\prime} = spec\, (E[X]/J)$ for two ideals $I \subset J$ of
$E[X].$ Put $A = E[X]/I$ and $B = E[X]/J$ and let \quad $\phi\,:\,
A \rightarrow B$ be the natural surjective ring homomorphism.
The homomorphism $\phi$ induces the surjective $E$-algebra homomorphism
$${\overline \phi}\,:\, {\overline A} \rightarrow {\overline B}$$
which upon taking fix points induces the $K$-algebra homomorphism
$${\overline \phi}^G\,:\, {\overline A}^G \rightarrow {\overline
B}^G. \tag{2.5}$$
By Remark 2.2  we see that ${\overline B}^G$ is
generated as a $K$-algebra by elements $\sum_{i=1}^n
\alpha_{j}^{\sigma_{i}} X^{\sigma_{i}}$ (more precisely their
images in ${{\overline B}}^G$). Similarly ${\overline A}^G$ is
generated as a $K$-algebra by elements $\sum_{i=1}^n
\alpha_{j}^{\sigma_{i}} X^{\sigma_{i}}$ (more precisely their
images in ${{\overline A}}^G$). It is clear that ${\overline
\phi}^G$ sends the element $\sum_{i=1}^n \alpha_{j}^{\sigma_{i}}
X^{\sigma_{i}} \in {{\overline A}}^G$ into $\sum_{i=1}^n
\alpha_{j}^{\sigma_{i}} X^{\sigma_{i}} \in {{\overline B}}^G.$
Hence ${\overline \phi}^G$ is onto. \qed\enddemo
\bigskip

\medskip
\proclaim{Lemma 2.6} Consider the group scheme $GL_r/E.$ The affine
group scheme $R_{E/K}(GL_{r})$ is a closed subscheme of
$GL_{rn}/K$
\endproclaim

\demo{\sl Proof} Let $\alpha_1,\dots,\alpha_n$ be a basis of $E$ over $K$ and let
Let $\beta_1,\dots,\beta_n$ be the corresponding dual basis with respect to $Tr_{E/K}.$
Define block matrices:

$$\A =
      \pmatrix  \alpha_{1}^{\sigma_{1}} I_r& \alpha_{1}^{\sigma_{2}} I_r& \dots & \alpha_{1}^{\sigma_{n}} I_r\\
                \alpha_{2}^{\sigma_{1}} I_r& \alpha_{2}^{\sigma_{2}} I_r& \dots & \alpha_{2}^{\sigma_{n}} I_r\\
                \vdots& \vdots& \dots&\vdots\\
                 \alpha_{n}^{\sigma_{1}} I_r& \alpha_{n}^{\sigma_{2}} I_r& \dots & \alpha_{n}^{\sigma_{n}} I_r
\endpmatrix ,
\quad
\BB  =
\pmatrix  \beta_{1}^{\sigma_{1}} I_r& \beta_{2}^{\sigma_{1}} I_r& \dots & \beta_{n}^{\sigma_{1}} I_r\\
                \beta_{1}^{\sigma_{2}} I_r& \beta_{2}^{\sigma_{2}} I_r& \dots & \beta_{n}^{\sigma_{2}} I_r\\
                \vdots& \vdots& \dots&\vdots\\
                 \beta_{1}^{\sigma_{n}} I_r& \beta_{2}^{\sigma_{n}} I_r& \dots & \beta_{n}^{\sigma_{n}} I_r
\endpmatrix
$$
Notice that by definition of the dual basis $\A \BB = \BB \A = I_{rn}.$
Define block diagonal matrices:
$$\X =
\pmatrix  X^{\sigma_{1}} & 0 I_r& \dots & 0 I_r\\
                0 I_r& X^{\sigma_{2}} & \dots & 0 I_r\\
                \vdots& \vdots& \dots&\vdots\\
                0 I_r& 0 I_r& \dots & X^{\sigma_{n}}
\endpmatrix ,
\quad
\Y =
\pmatrix  Y^{\sigma_{1}} & 0 I_r& \dots & 0 I_r\\
                0 I_r& Y^{\sigma_{2}} & \dots & 0 I_r\\
                \vdots& \vdots& \dots&\vdots\\
                0 I_r& 0 I_r& \dots & Y^{\sigma_{n}}
\endpmatrix ,$$
where $Y^{\sigma_{1}},\dots, Y^{\sigma_{n}}$ and $X^{\sigma_{1}},\dots, X^{\sigma_{n}},$
are multivariables written now in a form of $r \times r$ matrices indexed by $\sigma_{1},\dots, \sigma_{n}.$
Let $T_{ij}$ and $S_{ij},$ for all $1 \leq i \leq n, \, 1 \leq j \leq n,$ be $r \times r$
 multivariable matrices. Define block matrices of multivariables:
$$\T =
\pmatrix  T_{11}& T_{12}& \dots & T_{1n}\\
          T_{21}& T_{22}& \dots & T_{2n}\\
                \vdots& \vdots& \dots&\vdots\\
          T_{n1}& T_{n2}& \dots & T_{nn}
\endpmatrix ,
\quad
\SS =
\pmatrix  S_{11}& S_{12}& \dots & S_{1n}\\
          S_{21}& S_{22}& \dots & S_{2n}\\
                \vdots& \vdots& \dots&\vdots\\
          S_{n1}& S_{n2}& \dots & S_{nn}
\endpmatrix
$$
Notice that:
$$\A \X \BB =
      \pmatrix  \sum_{j=1}^n (\alpha_{1}\beta_{1})^{\sigma_{j}} X^{\sigma_{j}}&
                \sum_{j=1}^n (\alpha_{1}\beta_{2})^{\sigma_{j}} X^{\sigma_{j}}& \dots &
                \sum_{j=1}^n (\alpha_{1}\beta_{n})^{\sigma_{j}} X^{\sigma_{j}}\\
                \sum_{j=1}^n (\alpha_{2}\beta_{1})^{\sigma_{j}} X^{\sigma_{j}}&
                \sum_{j=1}^n (\alpha_{2}\beta_{2})^{\sigma_{j}} X^{\sigma_{j}}& \dots &
                \sum_{j=1}^n (\alpha_{2}\beta_{n})^{\sigma_{j}} X^{\sigma_{j}}\\
                \vdots& \vdots& \dots&\vdots\\
                \sum_{j=1}^n (\alpha_{n}\beta_{1})^{\sigma_{j}} X^{\sigma_{j}}&
                \sum_{j=1}^n (\alpha_{n}\beta_{2})^{\sigma_{j}} X^{\sigma_{j}}& \dots &
                \sum_{j=1}^n (\alpha_{n}\beta_{n})^{\sigma_{j}} X^{\sigma_{j}}
\endpmatrix $$
$$\A \Y \BB =
      \pmatrix  \sum_{j=1}^n (\alpha_{1}\beta_{1})^{\sigma_{j}} Y^{\sigma_{j}}&
                \sum_{j=1}^n (\alpha_{1}\beta_{2})^{\sigma_{j}} Y^{\sigma_{j}}& \dots &
                \sum_{j=1}^n (\alpha_{1}\beta_{n})^{\sigma_{j}} Y^{\sigma_{j}}\\
                \sum_{j=1}^n (\alpha_{2}\beta_{1})^{\sigma_{j}} Y^{\sigma_{j}}&
                \sum_{j=1}^n (\alpha_{2}\beta_{2})^{\sigma_{j}} Y^{\sigma_{j}}& \dots &
                \sum_{j=1}^n (\alpha_{2}\beta_{n})^{\sigma_{j}} Y^{\sigma_{j}}\\
                \vdots& \vdots& \dots&\vdots\\
                \sum_{j=1}^n (\alpha_{n}\beta_{1})^{\sigma_{j}} Y^{\sigma_{j}}&
                \sum_{j=1}^n (\alpha_{n}\beta_{2})^{\sigma_{j}} Y^{\sigma_{j}}& \dots &
                \sum_{j=1}^n (\alpha_{n}\beta_{n})^{\sigma_{j}} Y^{\sigma_{j}}
\endpmatrix .$$
Observe that the entries of $\A \X \BB$ and $\A \Y \BB$ are
$G$-equivariant. Hence, there is a well defined homomorphism of
$K$-algebras $$\Phi\,\,:\,\, K[\T, \SS]/(\T\SS - I_{rn},\, \SS\T -
I_{rn}) \,\,\, \rightarrow \,\,\, \bigl(M[\X, \Y]/(\X\Y -
I_{rn},\, \Y\X - I_{rn})\bigr)^G\tag{2.7}$$
$$\T \,\, \rightarrow\,\,  \A
\X \BB$$ $$\SS \,\, \rightarrow\,\,  \A \Y \BB$$
The definition of
$\Phi$ and the form of the entries of matrices $\A \X \BB$ and $\A
\Y \BB$ show (by the same argument as in Lemma 2.4) that the map
$\Phi$ is surjective. Observe that $$GL_{rn}/K = spec\,\, K[\T,
\SS]/(\T\SS - I_{rn},\, \SS\T - I_{rn}),$$ $$GL_r/E \,\, = \,\,
spec\,\, E[X, Y]/(XY - I_{r},\, YX - I_{r}),$$ where $X$ and $Y$ are
$r\times r$   multivariable matrices. Observe that there is a
natural $M$-algebra isomorphism $$M[\X, \Y]/(\X\Y - I_{rn},\, \Y\X
- I_{rn}) \cong
A^{\sigma_{1}}\otimes_{M}\dots\otimes_{M}A^{\sigma_{n}},$$ where
in this case $$A^{\sigma_{j}} = M[X, Y]/(XY - I_{r},\, YX -
I_{r})\,\,\cong\,\, M[X^{\sigma_{j}},
Y^{\sigma_{j}}]/(X^{\sigma_{j}}Y^{\sigma_{j}} - I_{r},\,
Y^{\sigma_{j}}X^{\sigma_{j}} - I_{r}).$$

\noindent
Hence, by Definition 2.3 we get a natural isomorphism of schemes over
$K:$ $$R_{E/K}(GL_r) \cong spec\, \bigl(M[\X, \Y]/(\X\Y -
I_{rn},\, \Y\X - I_{rn})\bigr)^G$$
and  it follows that $\Phi$ induces a closed imbedding of schemes over $K$
$$R_{E/K}(GL_r)\,\,\rightarrow GL_{rn}.\quad\qed$$\enddemo
\medskip

\noindent
\remark{\bf Remark 2.8} Let $E/K$ be an extension of two
local fields with the property that the corresponding extension of
rings of integers ${\cal O}_E/{\cal O}_K$ has integral basis
$\alpha_1,\dots,\alpha_n$ of ${\cal O}_E$ over ${\cal O}_K$ such
that the corresponding dual basis $\beta_1,\dots,\beta_n$ with
respect to $Tr_{E/K}$ is also a basis of ${\cal O}_E$ over ${\cal
O}_K.$ Let $R_{{\cal O}_E/{\cal O}_K}$ be the Weil restriction functor
defined analogously to the Weil restriction functor for the extension $E/K.$
Under these assumptions the following Lemmas 2.9 and 2.10 are
proven in precisely the same way as Lemmas 2.4 and 2.6.
\endremark
\medskip

\noindent \proclaim{Lemma 2.9} Let $V^{\prime} \subset V$ be a
closed imbedding of affine schemes over ${\cal O}_E.$ Then
$R_{{\cal O}_E/{\cal O}_K}(V^{\prime}) \subset R_{{\cal O}_E/{\cal
O}_K}(V)$ is a closed imbedding of affine schemes over
${\cal O}_K.$
\endproclaim

\medskip
\proclaim{Lemma 2.10} Consider the group scheme $GL_r/{\cal O}_E.$
The affine group scheme $R_{{\cal O}_E/{\cal O}_K}(GL_{r})$ is a
Zariski closed subscheme of $GL_{rn}/{\cal O}_K .$
\endproclaim
\medskip

We return to the case of the arbitrary separable field extension $E/K$ of degree $n.$
Every point of  $X_{0} \in GL_r(E)$ is uniquely determined  by the 
ring homomorphism
$$h_{X_{0}}\,:\, E[X, Y]/(XY - I_{r},\, YX - I_{r}) \,\,\rightarrow\,\, E$$
$$X  \,\mapsto\, X_0, \quad Y \,\mapsto\, Y_0,$$
where $Y_0$ is the inverse of $X_0.$
This gives immediately the homomorphism
$$h_{\T_{0}}\,:\,K[\T, \SS]/(\T\SS - I_{rn},\, \SS\T - I_{rn}) \rightarrow K$$
$$\T  \,\mapsto\, \T_{0} = \A \X_{0} \BB,$$
$$\SS \,\mapsto\, \SS_{0} = \A \Y_{0} \BB$$
where
$$\X_{0} =
\pmatrix  X_{0}^{\sigma_{1}} & 0 I_r& \dots & 0 I_r\\
                0 I_r& X_{0}^{\sigma_{2}} & \dots & 0 I_r\\
                \vdots& \vdots& \dots&\vdots\\
                0 I_r& 0 I_r& \dots & X_{0}^{\sigma_{n}}
\endpmatrix ,
\quad
\Y_{0} =
\pmatrix  Y_{0}^{\sigma_{1}} & 0 I_r& \dots & 0 I_r\\
                0 I_r& Y_{0}^{\sigma_{2}} & \dots & 0 I_r\\
                \vdots& \vdots& \dots&\vdots\\
                0 I_r& 0 I_r& \dots & Y_{0}^{\sigma_{n}}
\endpmatrix ,$$
and the action of $\sigma_i$ on $X_{0}$ and  $Y_{0}$ is the genuine
action on the entries of $X_{0}$ and  $Y_{0}.$ Obviously
$h_{\T_{0}}$ determines uniquely the point $\T_{0} \in GL_{rn}(K)$
with the inverse~$\SS_{0}.$
\medskip

\proclaim{Definition 2.11}
Assume that $Z = \{X_t;\,\, t \in T\} \subset GL_{r}(E)$ is a set of points.
We define the corresponding set of points:
$$Z_{\Phi}\quad = \quad \{\T_t = \A \X_{t} \BB; \,\, t \in T\} \quad \subset\quad  GL_{rn}(K).$$
We denote by $Z^{alg}$ the Zariski closure of $Z$ in $GL_{r}/E$ and by
$Z_{\Phi}^{alg}$ the Zariski closure of $Z_{\Phi}$ in
$GL_{rn}/K.$
\endproclaim

\medskip
\proclaim{Proposition 2.12} We have a natural isomorphism of schemes over~$K:$
$$R_{E/K}(Z^{alg}) \,\,\,\cong\,\,\, Z_{\Phi}^{alg}.$$
\endproclaim
\demo{\sl Proof} Let $$J_t = (XY - I_{r},\, YX - I_{r},\, X -
X_t,\, Y - Y_t)$$ be the prime ideal of $E[X,Y]$ corresponding to
the point $X_t \in GL_r(E).$ Let $$J = \bigcap_{t\in T} \, J_t.$$
By definition $Z^{alg} = spec\, (E[X,Y]/J).$ Let $$\J_{t} =
(\T\SS - I_{rn},\, \SS\T - I_{rn},\, \T - \A \X_{t} \BB,\, \SS -
\A \Y_{t} \BB )$$ be the prime ideal in $K[\T, \SS]/(\T\SS -
I_{rn},\, \SS\T - I_{rn})$ corresponding to the point $\A \X_{t}
\BB \in GL_{rn}(K).$ Define $$\J = \bigcap_{t\in T} \, \J_t.$$ By
definition $Z_{\Phi}^{alg} = spec\, (K[\T, \SS]/\J).$ Put $A = E[X,Y]/(XY - I_{r},\, YX - I_{r}).$ Observe that the
ring
${\overline A}^G$ is generated as a  K-algebra by $\A \X \BB$ and
$\A \Y \BB,$ since ${\overline A}$ is generated by $\X$ and $\Y$ as
an  $M$-algebra. Define $$\J_{t}^{\prime} = ( \A \X \BB - \A \X_{t}
\BB,\, \A \Y \BB  - \A \Y_{t} \BB )$$ which is an ideal of
${\overline A}^G.$ Put $$\J^{\prime} = \bigcap_{t\in T} \,
\J_{t}^{\prime}.$$ We have the following isomorphism induced by
$\Phi.$ $$K[\T, \SS]/\J_t \,\,\cong\,\, {\overline A}^G/\,\J_{t}^{\prime} \,\,\cong\,\, K\tag{2.13}$$
Hence, $\Phi^{-1}(\J_{t}^{\prime}) \,\,=\,\, \J_{t}$ and
$\Phi^{-1}(\J^{\prime}) \,\,=\,\, \J.$ This gives the isomorphism
$$K[\T, \SS]/\J \,\,\cong\,\, {\overline A}^G/\, \J^{\prime}.
\tag{2.14}$$
Let $B = E[X,Y]/J.$ There is a natural surjective
homomorphism of $K$-algebras coming from the construction in the
proof of Lemma 2.4 (see (2.5)): $${\overline A}^G/\, \J^{\prime}
\rightarrow {\overline B}^G \tag{2.15}$$
induced by the quotient map $A \rightarrow B.$ We want to prove that (2.15) is an isomorphism.
Observe that there is natural isomorphism of $K$-algebras:
$${\overline A}^G/\,\J_{t}^{\prime} \,\,\cong\,\, {\overline
{A/J_t}}^G \,\,\cong\,\, K. \tag{2.16}$$
Consider the following commutative diagram of homomorphisms of $K$-algebras: $$ \CD
{\overline A}^G/\,\J^{\prime} @>>> {\overline B}^G\\ @VVV @VVV\\
\prod_{t\in T} {\overline A}^G/\,\J_{t}^{\prime} @>{\cong}>>
\prod_{t\in T} {\overline {A/J_t}}^G
\endCD\tag{2.17}
$$ The left vertical arrow is an imbedding by definition of
$\J^{\prime}$ and the bottom horizontal arrow is an isomorphism by
(2.16). Hence the top horizontal arrow is an imbedding, i.e., the map
(2.15) is an isomorphism. The composition of maps (2.14) and (2.15) gives a natural
 isomorphism of $K$-algebras
$$K[\T, \SS]/\J \,\,\cong\,\,{\overline B}^G. \tag{2.18}$$
But $Z_{\Phi}^{alg} \, = \,  spec\,
(K[\T, \SS]/\J).$ In addition,  $Z^{alg}\,=\, spec\, B,$ hence
$R_{E/K}(Z^{alg}) \, = \, spec\, {\overline B}^G$ and Proposition 2.12
follows by (2.18). \qed\enddemo
\medskip

\remark{\bf Remark 2.19} If $Z$ is a subgroup of
$GL_{r}(E),$ then $Z_{\Phi}$ is a subgroup of $GL_{rn}(K).$ In
this case $Z^{alg}$ is a closed algebraic subgroup of
$GL_{r}/E$ and $Z_{\Phi}^{alg}$ is a closed algebraic subgroup
of $GL_{rn}/K.$
\endremark

\proclaim{Definition 2.20}
Let $G_1 = spec\, A$ be an affine algebraic group scheme
defined over $E$ and ${\frak g}_1$ its
Lie algebra. We define ${\frak g}=R_{E/K}{\frak g}_1$ to be the Lie algebra
obtained from ${\frak g}_1$ by considering it over $K$ with the same
bracket.
\endproclaim

\proclaim{Lemma 2.21} ${\Cal Lie}( R_{E/K}G_{1}) = R_{E/K}{\frak g}_1.$
\endproclaim
\demo{\sl Proof} Let $n=[E:K]$ and $G=Gal(E/K).$ Since $G_1$ is an
algebraic group ${\frak g}_1={\Cal Der}(A)$ is the Lie algebra of
derivations of the algebra $A$ of functions on
 $G_1$ [ H1].
Let $\phi : Der (A)\rightarrow Der (\bar A) $ be the  homomorphism
of Lie algebras (considered over $E$) given by the following formula: $${\phi}(\delta) =
\Sigma_{i=1}^{n} id\otimes\cdots\otimes id
\otimes\delta_{i}\otimes id\otimes \cdots\otimes id,$$ where
$\delta_{i}= \delta \otimes 1$ as an element of $Der(A^{\sigma
_i})$. Recall that $A^{\sigma _i}= A{\otimes}_{E,\sigma_i} M.$
 If $\sigma \in G$ and $\sigma (a_1\otimes\dots\otimes
a_n)=\sigma (a_{k_1})\otimes\dots\otimes\sigma(a_{k_n})$ one
readily sees that $\delta
_j(\sigma(a_{k_j}))=\sigma(\delta_{k_j})$ and therefore  $\phi
(\delta)$ is $G$-equivariant i.e., $\phi (\delta)\in Der ({\bar
A}^G).$ It is easy to see that $\phi(\delta)$ as an element of
$Der (\bar A)$ is nontrivial if $\delta$ is nontrivial. Since
$\phi(\delta)$ is $M$-linear and  ${\bar A}^G{\otimes}_{K} M =
\bar A,$ we see that $\phi(\delta)$ is a nontrivial element of $
Der ({\bar A}^G)={\Cal Lie }( R_{E/K}G_{1}).$ On the other hand,
observe that $${\Cal Lie}( R_{E/K}G_{1}){\otimes}_K{\bar K}= {\Cal
Lie}(R_{E/K}G_{1}{\otimes}_K{\bar K})=$$
$$ ={\Cal Lie}({\bar G}
\times_K \dots \times_K {\bar G})= (\oplus {\frak g}_1)\otimes_E
{\bar K}={\frak g}\otimes_K{\bar K}.$$
\medskip

\noindent
This  shows that ${\Cal Lie}( R_{E/K}G_{1})$ and $
R_{E/K}{\frak g}_1$ have the same dimensions and therefore are
equal.~ \qed\enddemo
\bigskip
\proclaim{Lemma  2.22} Let $\frak g$ be a Lie algebra over $E$ and let
${\frak g}^{\prime}$ be its derived algebra. Then
$$R_{E/K} ({\frak g}^{\prime}) = (R_{E/K} (\frak g))^{\prime}$$
\endproclaim
\demo{\sl Proof} This follows immediately from the fact that
$R_{E/K} (\frak g)$ and $\frak g$ have the same Lie bracket
(cf. Definition 2.20) \qed\enddemo

\proclaim{Lemma 2.23} If $G$ is an algebraic group over $E$ of characteristic 0, then
$$R_{E/K} (G^{\prime}) = (R_{E/K} G)^{\prime}$$
\endproclaim
\demo{\sl Proof} We have the following identities:

$$Lie ((R_{E/K} (G))^{\prime}) = (Lie (R_{E/K} (G)))^{\prime} =
(R_{E/K}(Lie(G)))^{\prime} = $$ $$ = R_{E/K}((Lie(G))^{\prime}) =
R_{E/K}(Lie(G^{\prime})) = Lie (R_{E/K}(G^{\prime}))$$ The first
and the fourth equality follow from Corollary on p.75 of [H1]. The
second and fifth equality follow from Lemma 2.21. The third
equality follows from Lemma 2.22. The Lemma follows by Theorem on
p. 87 of [H1] and Proposition on p. 110 of [H1]. \qed\enddemo

\subhead 3. Some remarks on bilinear forms
\endsubhead

\bigskip

\noindent Let $E$ be a finite extension of $\Q$ of degree $e.$ Let
$E_l = E\otimes \Q_l$ and  ${\cal O}_{E_l} = {\cal O}_{E} \otimes
\Z_l.$ Hence $E_l = \prod_{\lambda | l} E_{\lambda}$ and ${\cal
O}_{E_l} = \prod_{\lambda | l} O_{\lambda}.$ Let ${\cal
O}_{\lambda}^{\prime}$ be the dual to $ {\cal O}_{\lambda}$ with
respect to the trace $Tr_{E_{\lambda}/\Q_l}.$ For  $l \gg 0$ we have
${\cal O}_{\lambda}^{\prime} = {\cal O}_{\lambda}$ see [A], Chapter
7. From now on we take $l$ big enough to ensure that ${\cal
O}_{\lambda}^{\prime} = {\cal O}_{\lambda}$ and an abelian variety
$A$ has good reduction at all primes in ${\cal O}_E$ over $l.$ The
following lemma is the integral version of the sublemma 4.7 of [D].

\proclaim{Lemma 3.1} Let $T_1$ and $T_2$ be finitely generated, free
${\cal O}_{E_l}$-modules. For any bilinear $\Z_l$ (resp.
nondegenerate bilinear $\Z_l$) map
$$\psi_l\, :\, T_1 \times T_2 \rightarrow \Z_l$$
such that\,\, $\psi_l(ev_1, v_2) = \psi_l(v_1, ev_2)$\,\, for all $e
\in {\cal O}_{E_l}, v_1 \in T_1, v_2 \in T_2,$ there is a unique
${\cal O}_{E_l}$-bilinear (resp. nondegenerate ${\cal
O}_{E_l}$-bilinear ) map
$$\phi_l\, :\, T_1 \times T_2 \rightarrow {\cal O}_{E_l}$$
such that $Tr_{E_l/\Q_l}(\phi_l(v_1, v_2)) = \psi_l(v_1, v_2)$
for all $v_1 \in T_1$ and $v_2 \in T_2.$
\endproclaim

\demo{\sl Proof}
Similary to Sublemma 4.7, [D] we observe that the map
$$Tr_{E_l/\Q_l}\, :\, Hom_{{\cal O}_{E_l}} (T_1 \otimes_{{\cal O}_{E_l}}
T_2 \, ; {\cal O}_{E_l}) \rightarrow
Hom_{\Z_l} (T_1 \otimes_{{\cal O}_{E_l}} T_2 \, ; \Z_l)$$
is an isomorphism since it is a surjective map of torsion free $\Z_l$-modules
of the same $\Z_l$-rank.
The surjectivity of $Tr_{E_l/\Q_l}$ can be seen as follows.
The $\Z_l$-basis of the module $T_1{\otimes_{{\cal O}_{E_l}}} T_2$ is given by
$${\cal B} = \bigl( (0,\dots, 0, {\alpha}_{k}^{\lambda}, 0,\dots, 0)
e_i \otimes e^{\prime}_j \bigr)$$
where $(0,\dots, 0, {\alpha}_{k}^{\lambda}, 0,\dots, 0) \in
\prod_{\lambda | l} {\cal O}_{\lambda}$ and
${\alpha}_{k}^{\lambda}$ is an element
of a basis of ${\cal O}_{\lambda}$ over $\Z_l$ and
$e_i$ (resp. $e^{\prime}_j$) is an element of the standard basis of $T_1$
(resp. $T_2$) over ${\cal O}_{E_l}.$  Let
$\psi_{k,i,j}^{\lambda} \in  Hom_{\Z_l} (T_1 \otimes_{{\cal O}_{E_l}}
T_2 \, ;\, \Z_l)$ be the
homomorphism which takes value 1 on the element $(0,\dots, 0,
{\alpha}_{k}^{\lambda}, 0,\dots, 0) e_i \otimes e^{\prime}_j$
of the basis ${\cal B}$ and takes value 0 on the remaining elements of the basis
${\cal B}.$ Let us take
$\phi_{i,j} \in Hom_{{\cal O}_{E_l}} (T_1 \otimes_{{\cal O}_{E_l}} T_2 \,\, ; {\cal O}_{E_l})$
such that
$$\phi_{i,j} (e_r \otimes e_{s}^{\prime})  =
\cases  1 & \text{if $i = r$ and $j = s$ }\\
0   & \text{if $i \neq r$ or $j \neq s$}
\endcases
$$
Then for each $k$ there exist elements
(the dual basis) ${\beta}_{k}^{\lambda} \in {\cal O}_{\lambda}$
such that $Tr_{E_{\lambda}/\Q_l}({\beta}_{k}^{\lambda}
{\alpha}_{n}^{\lambda}) = \delta_{k,n}.$
If we put $\phi_{i,j,k}^{\lambda} = (0,\dots, 0, {\beta}_{k}^{\lambda}, 0,\dots, 0) \phi_{i,j}$
then clearly $Tr_{E_l/\Q_l}(\phi_{i,j,k}^{\lambda}(t_1, t_2)) =
\psi_{i,j,k}^{\lambda} (t_1, t_2).$ Hence the proof is finished since
the elements $\psi_{i,j,k}^{\lambda} (t_1, t_2)$ form a basis of
$Hom_{\Z_l} (T_1 \otimes_{{\cal O}_{E_l}} T_2 \, ; \Z_l)$  over $\Z_l.$
\qed\enddemo
\bigskip

Take $T_1 = T_2 = T_l$ and assume that $\psi_l$ is nondegenerate.
Let
$${\overline \psi_l}\, :\, T_l/l\,T_l \times T_l/l\,T_l \rightarrow \Z/l$$
be the $\Z/l$-bilinear pairing obtained by reducing the form $\psi_l$
modulo $l.$
Define
$$T_{\lambda} = e_{\lambda} T_l \cong T_l
\otimes_{{\cal O}_{E_l}} {\cal O}_{\lambda}$$
$$V_{\lambda} = T_{\lambda} \otimes_{{\cal O}_{\lambda}} E_{\lambda}$$
where $e_{\lambda}$ is the standard idempotent corresponding to
the decomposition ${\cal O}_{E_l} = \prod_{\lambda} {\cal O}_{\lambda}.$
Let $\pi_{\lambda}\, : \, {\cal O}_{E_l} \rightarrow  {\cal O}_{\lambda}$
be the projection.
We can define an ${\cal O}_{\lambda}$-nondegenerate bilinear form as follows:
$$\psi_{\lambda}\,:\, T_{\lambda} \times T_{\lambda} \rightarrow {\cal O}_{\lambda}$$
$$\psi_{\lambda}(e_{\lambda}v_1, e_{\lambda}v_2) = \pi_{\lambda} (\phi_l(v_1, v_2))$$
for any $v_1, v_2 \in T_l.$ Put $k_{\lambda} = {\cal O}_{\lambda}/\lambda\,{\cal O}_{\lambda}.$
This gives the $k_{\lambda}$-bilinear form
${\overline \psi}_{\lambda} = \psi_{\lambda} \otimes_{{\cal O}_{\lambda}} k_{\lambda}$
$${\overline \psi_{\lambda}}\, :\,
T_{\lambda}/\lambda\,T_{\lambda}\,\times\, T_{\lambda}/\lambda\,T_{\lambda}
\rightarrow k_{\lambda}.$$
We also have the
$E_{\lambda}$-bilinear form
$\psi_{\lambda}^0 := \psi_{\lambda} \otimes_{{\cal O}_{\lambda}} E_{\lambda}$
$$\psi_{\lambda}^0\,:\, V_{\lambda} \times V_{\lambda} \rightarrow
E_{\lambda}.$$

\proclaim{Lemma 3.2} Assume that the form ${\overline \psi}_l$ is
nondegenerate. Then the forms ${\overline \psi}_{\lambda},$
${\psi}_{\lambda}$ and ${\psi}_{\lambda}^0$ are nondegenerate for
each $\lambda | l.$
\endproclaim
\demo{\sl Proof} First we prove that ${\overline \psi}_{\lambda}$
is nondegenerate for all $\lambda | l.$
 Assume that ${\overline \psi}_{\lambda}$ is degenerate for some $\lambda.$
Without loss of generality we can assume that the left radical of
${\overline \psi}_{\lambda}$ is nonzero.
So there is a nonzero vector $e_{\lambda}v_0 \in T_{\lambda}$
(for some $v_0 \in T_l$) which maps to a nonzero
vector in $T_{\lambda}/\lambda\,T_{\lambda}$
such that $\psi_{\lambda}(e_{\lambda}v_0,\, e_{\lambda}w) \in \lambda\,{\cal O}_{\lambda}$
for all $w \in T_l.$
Now use the decomposition $T_l = \oplus_{\lambda} T_{\lambda},$\,\, Lemma 3.1
and the ${\cal O}_{E_l}$-linearity of $\phi_l$ to observe that for
each $w \in T_l$
$$\psi_l(e_{\lambda}v_0,\, w) =
Tr_{E_l/\Q_l}(\phi_l(e_{\lambda}v_0,\, \sum_{\lambda^{\prime}}
e_{\lambda^{\prime}} w)) = Tr_{E_{\lambda}/\Q_l} \psi_{\lambda}(e_{\lambda}v_0,\, e_{\lambda} w)
\in l \Z_l.$$
This contradicts the assumption that ${\overline \psi}_l$ is nondegenerate.

\noindent
Similarly, but in an easier way, we prove that ${\psi}_{\lambda}$ is nondegenerate.
From this immediately follows that ${\psi}_{\lambda}^0$ is nondegenerate.
\qed\enddemo

\subhead 4. Auxiliary facts about abelian varieties
\endsubhead

\noindent Let $A/F$ be a principally polarized, simple abelian variety of dimension $g$ with the
polarization defined over $F.$ Put ${\cal R} = End_{{\bar F}} (A)$ We assume that $End_{{\bar F}}
(A) = End_{F} (A)$,  hence the actions of ${\cal R}$ and $G_F$ on $A({\overline F})$ commute. Put $D
= End_{{\bar F}} (A) \otimes_{\Z} \Q.$ The ring ${\cal R}$ is an order in $D.$ Let $E_1$ be the
center of $D$ and let
$$E := \{a \in E_1;\,\, a^{\prime} = a\},$$
where ${\prime}$ is the Rosati involution. Let ${\cal R}_D$ be a maximal order in $D$ containing ${\cal R}.$
Put ${\cal O}_{E}^{0} := {\cal R} \cap E.$
The ring ${\cal O}_{E}^{0}$ is an order
in $E.$
Take $l$ that does not divide the index $[{\cal R}_D\, :\, {\cal R}].$ Then
${\cal R}_D \otimes_{\Z} \Z_l = {\cal R} \otimes_{\Z} \Z_l$ and
${\cal O}_{E} \otimes_{\Z} \Z_l = {\cal O}_{E}^{0} \otimes_{\Z} \Z_l$
\medskip

\noindent
The polarization of $A$ gives a $\Z_l$-bilinear,
nondegenerate, alternating pairing
$$\psi_l\,:\, T_l(A) \times T_l(A) \rightarrow \Z_l.\tag{4.1}$$
Because $A$ has the principal polarization, for any endomorphism
$\alpha \in {\cal R}$
we get ${\alpha}^{\prime} \in {\cal R},$ (see [Mi] chapter 13 and 17)
where ${\alpha}^{\prime}$ is the image of $\alpha$ by the Rosati involution.
Hence for any $v, w \in T_l(A)$ we have $\psi_l(\alpha v, w) = \psi_l(v, {\alpha}^{\prime}w)$
(see loc. cit.).

\medskip
\remark{\bf Remark 4.2}
When we deal with an
abelian variety which is not principally polarized, we must assume that
$l$ does not divide the degree of the polarization of $A,$ if we want
to get ${\alpha}^{\prime}\otimes 1 \in {\cal R}\otimes \Z_l$
for $\alpha \in {\cal R}.$
\endremark
\medskip

\noindent
By Lemma 3.1 there is a unique nondegenerate,
${\cal O}_{E_l}$-bilinear pairing
$$\phi_l\,:\, T_l(A) \times T_l(A) \rightarrow {\cal O}_{E_l} \tag{4.3}$$
such that $Tr_{E_l/\Q_l}(\phi_l(v_1, v_2)) = \psi_l(v_1, v_2).$
As in the general case define
$$T_{\lambda}(A) = e_{\lambda} T_l(A) \cong T_l(A)
\otimes_{{\cal O}_{E_l}} {\cal O}_{\lambda}$$
$$V_{\lambda}(A) = T_{\lambda}(A) \otimes_{{\cal O}_{\lambda}}
{E}_{\lambda}.$$
Note that $T_{\lambda}(A)/ \lambda T_{\lambda}(A) \cong A[\lambda]$
as $k_{\lambda}[G_F]$-modules.

\noindent
Again as in the general case define nondegenerate, ${\cal O}_{\lambda}$-bilinear form
$$\psi_{\lambda}\,:\, T_{\lambda}(A) \times T_{\lambda}(A) \rightarrow
{\cal O}_{\lambda} \tag{4.4}$$
$$\psi_{\lambda}(e_{\lambda}v_1, e_{\lambda}v_2) =
\pi_{\lambda} (\phi_l(v_1, v_2))$$
for any $v_1, v_2 \in T_l(A).$ The form $\psi_{\lambda}$ gives the forms:
$${\overline \psi}_{\lambda}\,:\, A[\lambda] \times A[\lambda] \rightarrow
k_{\lambda}. \tag{4.5}$$
$${\psi}_{\lambda}^0\,:\, V_{\lambda}(A) \times V_{\lambda}(A) \rightarrow
E_{\lambda}. \tag{4.6}$$

\noindent
Notice that all bilinear forms $\psi_{\lambda},
{\overline \psi}_{\lambda}$ and
$\psi_{\lambda}^0$ are alternating forms. For $l$ relatively prime to the degree
of polarization the form ${\psi}_{l}$ is nondegenerate. Hence by lemma 3.2 the forms
$\psi_{\lambda}, {\overline \psi}_{\lambda}$ and $\psi_{\lambda}^0$ are nondegenerate.

\proclaim{Lemma 4.7} Let $\chi_{\lambda}\,:\, G_F \rightarrow \Z_l
\subset {\cal O}_{\lambda}$
be the composition of the cyclotomic character with the natural imbedding
$\Z_l \subset {\cal O}_{\lambda}.$

\roster
\item"{(i)}"\quad For any $\sigma \in G_F$ and all $v_1, v_2 \in T_{\lambda}(A)$
$$\psi_{\lambda}(\sigma v_1, \sigma v_2) = \chi_{\lambda}(\sigma) \psi_{\lambda}( v_1, v_2). $$
\item"{(ii)}"\quad For any $\alpha \in {\cal R}$ and all $v_1, v_2 \in T_{\lambda}(A)$
$$\psi_{\lambda}(\alpha v_1,  v_2) =  \psi_{\lambda}( v_1, {\alpha}^{\prime} v_2). $$
\endroster
\endproclaim
\demo{\sl Proof} Proof is  the same as the proof of Lemma 2.3 in [C2].
\qed\enddemo

\medskip
\remark{\bf Remark 4.8}
After tensoring appropriate objects with $\Q_l$ in lemmas 3.1 and 4.6
we obtain Lemmas 2.2 and 2.3 of [C2].
\endremark
\medskip

Up to theorem 4.16, $A/F$ is an arbitrary
abelian variety defined over a number field $F.$ We introduce
some notation. Let $G_{l^{\infty}},$ $G_l,$ $G_{l^{\infty}}^0$ denote
the images of the corresponding representations:
$$\rho_l\, :\, G_F \rightarrow GL(T_l(A)) \cong GL_{2g}(\Z_l),$$
$${\overline{\rho_l}}\, :\, G_F \rightarrow GL(A[l]) \cong
GL_{2g}(\F_l),$$
$$\rho_l \otimes \Q_l \, :\, G_F \rightarrow
GL(V_l(A)) \cong GL_{2g}(\Q_l).$$
\medskip

\noindent
Let  ${\cal G}_{l}^{alg},$ (\, $G(l)^{alg},$\,\,\,  $G_{l}^{alg}$\,\,\, resp.) denote
the Zariski closure of the image of the representation $\rho_l,$ (\, ${\overline{\rho_l}},$\,\,\,
$\rho_l \otimes \Q_l$\,\,\, resp.) in $GL_{2g}/\Z_{l},$ (\, $GL_{2g}/\F_{l},$\,\,\, $GL_{2g}/\Q_{l}$
\,\,\, resp.)
\medskip

\noindent
Due to our assumptions on the $G_F$-action and the properties of
the forms $\psi_{\lambda}, {\overline
\psi}_{\lambda}$ and $\psi_{\lambda}^0$ we get:

$$G_{l^{\infty}}  \,\,\subset\,\,
{\cal G}_{l}^{alg}(\Z_l) \,\,\subset\,\, \prod_{\lambda | l} \,\,
GSp_{T_{\lambda}(A)}({\cal O}_{\lambda}) \,\,\subset\,\,
GL_{T_l(A)}(\Z_l) \tag{4.9}$$
$$G_l \,\,\subset\,\, G(l)^{alg}(\F_l)
\,\,\subset \,\, \prod_{\lambda | l} \,\,
GSp_{A[\lambda]}(k_{\lambda}) \,\,\subset  \,\, GL_{A[l]}(\F_l) \tag{4.10}$$
$$G_{l^{\infty}}^0 \,\,\subset\,\, G_{l}^{alg}(\Q_l)
\,\,\subset\,\, \prod_{\lambda | l} \,\,
GSp_{V_{\lambda}(A)}(E_{\lambda}) \,\,\subset\,\,
GL_{V_l(A)}(\Q_l). \tag{4.11}$$
\bigskip

\noindent
Let $K/\Q_l$ be a local field extension and ${\cal O}_K$ the
ring of integers in $K.$
Let $T$ be a finitely generated, free ${\cal O}_K$-module and
let $V = T \otimes_{{\cal O}_K} K.$
Consider a continuous representation
$\rho\,:\, G_F\, \rightarrow GL(T)$ and the induced
representation $\rho^0\,=\, \rho \otimes K\,\, :\,\, G_F \rightarrow GL(V).$
Since $G_F$ is compact and $\rho^0$ is continuous, the subgroup $\rho^0 (G_F)$
of $GL(V)$ is closed. By [Se7], LG. 4.5, $\rho^0 (G_F)$ is an analytic
subgroup of $GL(V).$

\proclaim{Lemma 4.12} Let ${\frak g}$ be the Lie algebra of
the group $\rho^0 (G_F)$
\roster
\item"{(i)}"\quad There is an open subgroup $U_0 \subset \rho^0(G_F)$ such that
$$End_{U_0}\, (V) = End_{\frak g}\, (V).$$
\item"{(ii)}"\quad For all open subgroups $U \subset \rho^0(G_F)$ we have
$$End_{U}\, (V) \subset End_{\frak g}\, (V).$$
\item"{(iii)}"\quad Taking union over all open subgroups $U \subset \rho^0(G_F)$ we get
$$\bigcup_{U}\, End_{U}\, (V) = End_{\frak g}\, (V).$$
\endroster
\endproclaim
\demo{\sl Proof} (i)\,\,\,
Note that for any open subgroup ${\tilde U}$ of  ${\frak g}$ we have
$$End_{{\tilde U}}\, (V) \, =\,  End_{\frak g}\, (V) \tag{4.13}$$
because $K\, {\tilde U}\, = \, {\frak g}.$ By [B], Prop. 3, III.7.2,
for some open ${\tilde U} \subset {\frak g},$ there is an exponential map
$$exp\, :\, {\tilde U}\, \rightarrow \rho^0(G_F)$$
which is an analytic isomorphism and such that $exp\, ({\tilde U})$ is an
open subgroup
of $\rho^0(G_F).$ The exponential map can be expressed by
the classical power
series for $exp\, (t).$ On the other hand
by [B], Prop. 10, III.7.6, for some open
$U \subset \rho^0 (G_F),$ there is a logarithmic map
$$log\, :\,  U\, \rightarrow {\frak g}$$
which is an analytic isomorphism and the inverse of $exp.$ The logarithmic
map can be expressed by the classical power series for $\ln t.$ Hence, we
can choose ${\tilde U}_0$ such that
$U_0 = exp\, ({\tilde U}_0)$ and $log\, (U_0) = {\tilde U}_0.$ This gives
$$End_{U_{0}}\, (V) \, =\, End_{{\tilde U}_{0}}\, (V). \tag{4.14}$$
and (i) follows by (4.13) and (4.14).
\medskip

\noindent
(ii)\,\,\, Observe that for any open $U \subset \rho^0 (G_F)$ we have
$$End_{U}\, (V) \, \subset \, End_{ U_0 \cap U}\, (V).$$
Hence \, (ii)\, follows by \, (i).
\medskip

\noindent
(iii)\,\,\, Follows by \, (i) \, and \, (ii).
\qed\enddemo
\bigskip

\noindent
\proclaim{Lemma 4.15} Let $A/F$ be an abelian variety over $F$ such that
$End_{F}\, (A) = End_{{\overline F}}\, (A).$
Then
$$End_{{G_F}}\, (V_l(A))= End_{\frak g_l}\, (V_l(A)).$$
\endproclaim
\demo{\sl Proof} By the result of Faltings [Fa], Satz 4,
$$End_{{L}}\, (A) = End_{G_L}\, (V_l(A))$$
for any finite extension $L/F.$ By the assumption
$End_{F}\, (A) = End_{L}\, (A).$ Hence
$$End_{{G_F}}\, (V_l(A))= End_{U^{\prime}}\, (V_l(A))$$
for any open subgroup $U^{\prime}$ of $G_F.$
So the claim follows by Lemma 4.12 (iii).
\qed\enddemo
\bigskip

Let $A$ be a simple abelian variety defined over $F$ and $E$ be the center of the algebra $D=End_F(A)\otimes {\Q}.$
Let ${\lambda}|l$ be a prime of ${\cal O}_E$ over $l.$ Consider the following representations.
\medskip

\noindent
$$\rho_{\lambda}\, :\, G_F \rightarrow GL(T_{\lambda}(A)),$$
$${\overline{\rho_{\lambda}}}\, :\, G_F \rightarrow GL(A[\lambda]),$$
$$\rho_{\lambda} \otimes_{{\cal O}_{\lambda}} E_{\lambda} \, :\, G_F
\rightarrow GL(V_{\lambda}(A)),$$
\medskip

\noindent
where $\lambda |l.$ Let ${\cal G}_{\lambda}^{alg},$ (\, $G(\lambda)^{alg},$\,\,\,  $G_{\lambda}^{alg}$\,\,\, resp.) denote
the Zariski closure of the image of the representation $\rho_{\lambda},$ (\,
${\overline{\rho_\lambda}},$\,\,\,
$\rho_{\lambda} \otimes E_{\lambda}$\,\,\, resp.) in $GL_{T_{\lambda}(A)}/{\cal O}_{\lambda},$
(\, $GL_{A[\lambda]}/k_{\lambda},$\,\,\, $GL_{V_{\lambda}(A)}/E_{\lambda}$
\,\,\, resp.)
\medskip

\noindent
\proclaim{Theorem 4.16} Let A be a simple abelian variety with the property that
${\cal R} = End_{{\bar F}} (A) = End_{F} (A).$
Let ${\cal R}_{\lambda} = {\cal R} \otimes_{{\cal O}_E^{0}} {\cal O}_{\lambda}$
and let $D_{\lambda} = D \otimes_{E} E_{\lambda}.$ Then
\roster
\item"{(i)}"\quad
$End_{{\cal O}_{\lambda}[G_F]}\, (T_{\lambda}(A)) \, \cong\, {\cal R}_{\lambda} $
\item"{(ii)}"\quad
$End_{R_{\lambda}[G_F]}\, (V_{\lambda}(A)) \, \cong\, D_{\lambda}$
\item"{(iii)}"\quad
$End_{k_{\lambda}[G_F]}\, (A[\lambda]) \, \cong\, {\cal R}_{\lambda} \otimes_{{\cal O}_{\lambda}}
k_{\lambda}$
\endroster
\endproclaim
\demo{\sl Proof} It follows by [Fa], Satz 4 and [Za],
Cor. 5.4.5.
\qed\enddemo

\proclaim{Lemma 4.17} Let $K$ be a field and let
$R$ be a unital $K$-algebra. Put $D = End_{R} (M)$ and let
$L$ be a subfield of the center of $D.$ Assume that $L/K$ is a
finite separable extension. If $M$ is a semisimple $R$-module
then $M$ is also a semisimple $R \otimes_{K} L$-module with
the obvious action of $R \otimes_{K} L$ on $M.$
\endproclaim
\demo{\sl Proof} Take $\alpha \in L$ such that $L = K(\alpha).$
Let $[L: K] = n.$ Let us write $M = \oplus_{i} M_i$ where all $M_i$ are
simple $R$ modules. For any $i$ we put ${\tilde M}_i =
\sum_{k=0}^{n-1} \alpha^k\, M_i.$ Then ${\tilde M}_i$ is an
$R \otimes_{K} L$-module.
Because $M_i$ is a simple $R$-module we can write
$${\tilde M}_i = \bigoplus_{k=0}^{m-1} \alpha^k\, M_i,$$
for some $m.$ Observe that if $m=1,$ then ${\tilde M}_i$
is obviously a simple $R \otimes_{K} L$-module. If $m{>}1,$ we pick
any simple $R$-submodule $N_i \subset {\tilde M}_i,$ $N_i \not= M_i.$
There is an $R$- isomorphism
$\phi \,:\, M_i \rightarrow N_i$ by  semisimplicity of ${\tilde M}_i.$
 We can write $M = M_i \oplus N_i \oplus M^{\prime},$
where $M^{\prime}$ is an $R$-submodule of $M.$
Define $\Psi \in Aut_{R} (M) \subset End_{R} (M)$ by
$\Psi|_{M_i} = \phi,$ $\Psi|_{N_i} = \phi^{-1}$ and
$\Psi|_{M^{\prime}} = Id_{M^{\prime }}.$ Note that
$$\Psi (\bigoplus_{k=0}^{m-1} \alpha^k\, M_i) =
\bigoplus_{k=0}^{m-1} \alpha^k\, N_i \tag{4.18}$$
since $\alpha$ is in the center
of $D.$ Hence ${\tilde M_i} =
\bigoplus_{k=0}^{m-1} \alpha^k\, N_i$
by the classification of semisimple $R$-modules.
We conclude that ${\tilde M}_i$ is a simple $R \otimes_{K} L$-module.
Indeed, if $N \subset {\tilde M}_i$
were a nonzero $R \otimes_{K} L$-submodule of
${\tilde M}_i$ then we could pick any simple $R$-submodule
$N_i \subset N.$ If $N_i = M_i$ then $N = {\tilde M}_i.$ If
$N_i \not= M_i$ then by (4.18)
${\tilde M}_i = \bigoplus_{k=0}^{m-1} \alpha^k\, N_i \subset N.$
Since $M = \sum_{i} {\tilde M}_i,$ we see that $M$ is a
semisimple $R \otimes_{K} L$-module.
\qed\enddemo

\proclaim{Theorem 4.19}  Let A be a simple abelian variety with the property that ${\cal R} =
End_{{\bar F}} (A) = End_{F} (A).$ Let ${\cal R}_{\lambda} = {\cal R} \otimes_{{\cal O}_E^{0}}
{\cal O}_{\lambda}$ and let $D_{\lambda} = D \otimes_{E} E_{\lambda}.$ Then $G_F$ acts on
$V_{\lambda}(A)$ and $A[\lambda]$ semisimply and $G_{\lambda}^{alg}$ and $G(\lambda)^{alg}$ are
reductive algebraic groups. The scheme ${\cal G}_{\lambda}^{alg}$ is a reductive group scheme over
${\cal O}_{\lambda}$ for $l$ big enough.
\endproclaim
\demo{\sl Proof} It follows by [Fa], Theorem 3 and our Lemma 4.17. The last statement follows by
[LP1], Proposition 1.3, see also [Wi], Theoreme 1. \qed\enddemo
\medskip

\subhead 5. Abelian varieties of type I and II
\endsubhead

\noindent In this section we work with abelian varieties of type I
and II in the Albert's classification list of  division algebras
with involution [M], p. 201, i.e. $E \subset D = End_{\bar F} (A)
\otimes_{Z} \Q$ is the center of $D$ and $E$ is a totally real
extension of $\Q$ of degree $e.$ Moreover $D$ is either $E$ or an
indefinite quaternion algebra with the center $E,$ such that $D
\otimes_{E} \R \,\cong  \,M_{2,2}(\R).$ In the first part of this
section we prove integral versions of the results of Chi [C2] for
abelian varieties of type II. Let $l$ be any prime number that
does not divide the index $[{\cal R}_D\, :\, {\cal R}]$ and $l \gg
0$ for which $D \otimes_{\Q} \Q_l$ splits over $\Q_l.$ With such a
choice of $l,$ the algebra $D$ splits over any prime $\lambda$ in
${\cal O}_E$ over $l$. Hence, $D_{\lambda}=M_{2,2}(E_{\lambda}).$
Then by [R, Corollary 11.2 p. 132 and Theorem 11.5 p. 133] the
ring $R_{\lambda}$ is a maximal order in $D_{\lambda}.$ So by [R]
Theorem 8.7 p. 110 we get $R_{\lambda}=M_{2,2}({\Cal
O}_{\lambda})$, hence $R_{\lambda}{\otimes}_ {{\Cal O}_{\lambda}}
k_{\lambda}=M_{2,2}(k_{\lambda}).$ Similarly to [C2] we put
$$t= \pmatrix  1&0\\
    0&-1
 \endpmatrix ,
\quad
u= \pmatrix  0&1\\
    1&0
 \endpmatrix .
$$
Let $e=\frac{1}{2}(1+t),$ $f=\frac{1}{2}(1+u),$ \,\, ${\Cal X}=e\, T_{\lambda}(A),$
${\Cal Y}=(1-e)\, T_{\lambda}(A),$ ${\Cal X^{\prime}}=f\, T_{\lambda}(A),$
${\Cal Y^{\prime}}=(1-f)\, T_{\lambda}(A).$
Put $X={\Cal X}{\otimes}_{{\Cal O}_{\lambda}} E_{\lambda},$
$X^{\prime}={\Cal X}^{\prime}{\otimes}_{{\Cal O}_{\lambda}} E_{\lambda},$
 $Y={\Cal Y}{\otimes}_{{\Cal O}_{\lambda}} E_{\lambda},$
$Y^{\prime}={\Cal Y}^{\prime}{\otimes}_{{\Cal O}_{\lambda}} E_{\lambda},$
$\overline{\Cal X}={\Cal X}{\otimes}_{{\Cal O}_{\lambda}} k_{\lambda},$
$\overline{\Cal X}^{\prime}={\Cal X}^{\prime}{\otimes}_{{\Cal O}_{\lambda}} k_{\lambda},$
 $\overline{\Cal Y}={\Cal Y}{\otimes}_{{\Cal O}_{\lambda}} k_{\lambda},$
$\overline{\Cal Y}^{\prime}={\Cal Y}^{\prime}{\otimes}_{{\Cal O}_{\lambda}} k_{\lambda}.$
Multiplication by $u$ gives an
${\Cal O}_{\lambda}[G_F]$-isomorphism between ${\Cal X}$ and ${\Cal Y},$
hence it yields an $E_{\lambda}[G_F]$-isomorphism between $X$ and $Y$ and
a $k_{\lambda}[G_F]$-isomorphism between $\overline{\Cal X}$ and
$\overline{\Cal Y}.$ Multiplication by $t$ gives
an ${\Cal O}_{\lambda}[G_F]$-isomorphism between ${\Cal X}^{\prime}$ and
${\Cal Y}^{\prime},$
hence it yields an $E_{\lambda}[G_F]$-isomorphism between $X^{\prime}$ and
$Y^{\prime}$ and a
$k_{\lambda}[G_F]$-isomorphism between $\overline{\Cal X}^{\prime}$ and
$\overline{\Cal Y}^{\prime}.$ Observe that
$$End_{{\cal O}_{\lambda}[G_F]}\, ({\Cal X}) =
End_{{\cal O}_{\lambda}[G_F]}({\Cal X}^{\prime}) = {\Cal O}_{\lambda} \tag{5.1}$$
$$End_{E_{\lambda}[G_F]}\, (X) =
End_{E_{\lambda}[G_F]}(X^{\prime}) = E_{\lambda}\tag{5.2}$$
$$End_{k_{\lambda}[G_F]}\, (\overline{\Cal X}) =
End_{k_{\lambda}[G_F]}(\overline{\Cal X}^{\prime}) = k_{\lambda}. \tag{5.3}$$

\noindent
So the representations of $G_F$ on the spaces $X , Y , X^{\prime},
Y^{\prime}$ (resp. $\overline{\Cal X} ,\overline{\Cal Y} ,
\overline{\Cal X}^{\prime},
\overline{\Cal Y}^{\prime}$) are absolutely irreducible over
$E_{\lambda}$ (resp. over $k_{\lambda}$). Hence, the bilinear form
$\psi_{\lambda}^{0}$ cf. (4.4) (resp. ${\overline \psi}_{\lambda}$ cf. (4.5))
when restricted to any of the spaces $X, X^{\prime}, Y,
Y^{\prime},$ (resp. spaces $\overline{\Cal X}, {\overline{\Cal X}}^{\prime},
\overline{\Cal Y}, {\overline{\Cal Y}}^{\prime}$)
is either nondegenerate or isotropic.
\medskip

\noindent We obtain the integral version of Theorem A of [C2].

\proclaim{Theorem 5.4} If $A$ is of type II, then there is a free
${\cal O}_{\lambda}$-module ${\cal W}_{\lambda}(A)$ of rank $2h$ such
that
\roster
\item"{(i)}" we have an isomorphism of ${\cal O}_{\lambda}[G_{F}]$-
modules $T_{\lambda}(A) \cong {\cal W}_{\lambda}(A)
\oplus {\cal W}_{\lambda}(A)$
\item"{(ii)}" there is an alternating pairing
$\psi_{\lambda}\,:\, {\cal W}_{\lambda}(A) \times {\cal W}_{\lambda}(A)
\rightarrow {\cal O}_{\lambda}$
\item"{(ii')}" the induced alternating pairing
$\psi_{\lambda}^{0} \,:\,  W_{\lambda}(A) \times W_{\lambda}(A) \rightarrow E_{\lambda}$ is
nondegenerate, where  $W_{\lambda}(A) = {\cal W}_{\lambda}(A) \otimes_{{\cal
O}_{\lambda}} E_{\lambda}$
\item"{(ii'')}" the induced alternating pairing
${\overline \psi_{\lambda}}\,:\,  {\overline {\cal W}_{\lambda}}(A) \times {\overline {\cal
W}_{\lambda}}(A) \rightarrow k_{\lambda}$ is nondegenerate, where ${\overline {\cal
W}_{\lambda}}(A) = {\cal W}_{\lambda}(A) \otimes_{{\cal O}_{\lambda}} k_{\lambda}.$
\endroster
The pairings in (ii), (ii') and (ii'') are compatible with  the $G_F$-action in the
same way as the pairing in Lemma 4.7 (i).
\endproclaim

\demo{\sl Proof} (ii') is proven in [C2], (i) and (ii) are
straightforward generalizations  of the arguments in {\sl loc. cit}.
The bilinear pairing ${\phi}_{l}$ is nondegenerate, hence the bilinear pairing
${\overline \phi}_{l}$ is nondegenerate, since the abelian variety
$A$ is principally polarized by assumption.
Actually ${\overline \phi}_{l}$ is
nondegenerate for any abelian variety with polarization degree
prime to $l.$ So, by Lemma 3.2 the form ${\overline \psi}_{\lambda}$
is nondegenerate for all $\lambda$ hence simultaneously the forms
$\psi_{\lambda}^{0}$ and ${\overline \psi_{\lambda}}$ are
 nondegenerate. Now we finish the proof of (ii'')
arguing  for  $A[\lambda]$ similarly as it is done for
$V_{\lambda}$ in [C2], Lemma 3.3.
\qed\enddemo

\bigskip
From now on we work with the abelian varieties of type either I or II.
We assume that the field $F$ of definition of $A$ is such that
$G_{l}^{alg}$ is a connected algebraic group.

\noindent
Let us put
$$ T_{\lambda}  = \cases T_{\lambda}(A)   & \text{if $A$ is of type
I }\\
\\{\Cal W}_{\lambda}(A)\, ,  & \text{if $A$ is of type
II }
\endcases \tag{5.5}$$
Let  $V_{\lambda}=T_{\lambda}\otimes_{{\cal O}_{\lambda}} E_{\lambda}$ and
$A_{\lambda}\, =\, V_{\lambda} / T_{\lambda}.$
 With this notation we have:
$$ V_l(A)  = \cases \bigoplus_{\lambda | l} V_{\lambda}   & \text{if $A$ is of type
I }\\
\\ \bigoplus_{\lambda | l} \bigl(V_{\lambda} \oplus V_{\lambda}\bigr)\, ,  & \text{if $A$ is of type
II }
\endcases \tag{5.6}$$
We put
$$V_l \,\,=\,\, \bigoplus_{\lambda | l} V_{\lambda}\tag{5.7}$$
Let $V_{\Phi_{\lambda}}$ be the space $V_{\lambda}$ considered over
$\Q_l.$ For $\rho_{\lambda}(g) = X_{\lambda}$ and for $ X_{\lambda}
\in GL(V_{\lambda}),$ we define $\rho_{\Phi_{\lambda}}(g) \, = \,
\T_{\lambda}\, = \, \A_{\lambda}\X_{\lambda}\BB_{\lambda}$ ( cf. the
definition of the map $\Phi$ in (2.7) for the choice of
$\A_{\lambda}$ and $\BB_{\lambda} ).$ We have the following equality
of $\Q_l$-vector spaces:
$$V_l \,\,=\,\, \bigoplus_{\lambda | l} V_{\Phi_{\lambda}}\tag{5.8}$$
The $l$-adic representation
$$\rho_{l}\, :\, G_F \longrightarrow GL(V_l(A)) \tag{5.9}$$
induces the following representations (note that we use the notation $\rho_{l}$ for both representations
  (5.9) and (5.10) cf. Remark 5.13 ):
$$\rho_{l}\,\colon\,G_F\longrightarrow GL(V_l) \tag{5.10}$$
$$
\prod{\rho}_{\lambda}\,\colon\,G_F\longrightarrow \prod_{\lambda}GL(V_{\lambda})
 \tag{5.11}$$
$$
 {\prod{\rho}_{\Phi_{\lambda}}}\,\colon\,G_F\longrightarrow {\prod}_{\Phi_{\lambda}}
 GL(V_{\Phi_{\lambda}}).
\tag{5.12}$$
\medskip

\noindent \remark{\bf Remark 5.13} In the case of abelian variety of type II we have $V_l(A) =
V_{l} \oplus V_{l}$ and the action of $G_F$ on the direct sum is the diagonal one as follows from
Theorem 5.4. Hence, the images of the Galois group via the representations (5.9), (5.10) and
(5.12) are isomorphic. Also the Zariski closures of the images of these three representations are
isomorphic as algebraic varieties over $\Q_l$ in the corresponding $GL$-groups. Similarly,
$V_{\lambda}(A) = V_{\lambda} \oplus V_{\lambda}$ with the diagonal action of $G_F$ on the direct
sum by Theorem 5.4. Hence, the images of the representations given by $G_F$-actions on
$V_{\lambda}$ and $V_{\lambda}(A)$ are isomorphic and so are their Zariski closures in
corresponding $GL$-groups. For this reason, in the sequel, we will identify the representation of
$G_F$ on $V_l(A)$ (respectively on $V_{\lambda}(A)$) with its representation on $V_l$ (resp.
$V_{\lambda}$).
\endremark
\bigskip

\noindent By Remark 5.13 we can consider $G_{l}^{alg}$\, \, (resp. $G_{{\lambda}}^{alg}$) \, to be
the Zariski closure in $GL_{V_{l}}$ (resp. $GL_{V_{{\lambda}}}$) of the image of the
representation $\rho_{l}$ of (5.10) (resp. $\rho_{\lambda}$ of (5.11)). Let \,
$G_{\Phi_{\lambda}}^{alg}$\, denote the Zariski closure in $GL_{V_{\Phi_{\lambda}}}$ of the image
of the representation ${\rho}_{\Phi_{\lambda}}$ of (5.12). Let ${\frak g}_{l}$ be the Lie algebra
of $G_{l}^{alg},$ ${\frak g}_{\lambda}$ be the Lie algebra of $G_{\lambda}^{alg}$ and let ${\frak
g}_{\Phi_{\lambda}}$ be the Lie algebra of $G_{\Phi_{\lambda}}^{alg}.$  By definition, we have the
following inclusions:
$$G_{l}^{alg}\,\subset\,
{\prod}_{{\lambda}|l}G_{{\Phi}_{\lambda}}^{alg} \tag{5.14}$$
$$(G_{l}^{alg})^{\prime}\,\subset\,
{\prod}_{{\lambda}|l}(G_{{\Phi}_{\lambda}}^{alg})^{\prime}
\tag{5.15}$$
$${\frak g}_{l}\,\subset\,
{\bigoplus}_{{\lambda}|l}{\frak g}_{{\Phi}_{\lambda}} \tag{5.16}$$
$${\frak g}_{l}^{ss}\,\subset\,
{\bigoplus}_{{\lambda}|l}{\frak g}_{{\Phi}_{\lambda}}^{ss}. \tag{5.17}$$

\noindent
The map (5.14) gives a map
$$G_{l}^{alg} \rightarrow  G_{{\Phi}_{\lambda}}^{alg},\tag{5.18}$$
which induces the natural map of Lie algebras:
$${\frak g}_{l} \rightarrow {\frak g}_{{\Phi}_{\lambda}}.\tag{5.19}$$

\proclaim{Lemma 5.20} The map (5.19)  of Lie algebras is surjective for any prime $\lambda | l.$
Hence the following map of Lie algebras:
$${\frak g}_{l}^{ss} \rightarrow {\frak g}_{{\Phi}_{\lambda}}^{ss} \tag{5.21}$$
is surjective.
\endproclaim
\demo{\sl Proof} We know by the result of Tate, [T2] that the $\Q_l[G_F]$-module $V_l(A)$ is of
the Hodge-Tate type for any prime $v$ of ${\cal O}_F$ dividing $l.$ Hence by the theorem of
Bogomolov cf. [Bo] we have
$${\frak g}_{l} = {\cal Lie}\, (\rho_{l}(G_F)).$$
Since each $\Q_l[G_F]$-module $V_{{\Phi}_{{\lambda}}}$ is a direct summand of the
$\Q_l[G_F]$-module $V_l,$ then the $\Q_l[G_F]$-module $V_{{\Phi}_{{\lambda}}}$ is also of the
Hodge-Tate type for any prime $v$ of ${\cal O}_F$ dividing $l.$ It follows by the theorem of
Bogomolov, [Bo] that
$${\frak g}_{{\Phi}_{\lambda}} = {\cal Lie}\, (\rho_{{\Phi}_{\lambda}}(G_F)).$$
But the surjective map of $l$-adic Lie groups $\rho_{l}(G_F)
\rightarrow \rho_{{\Phi}_{\lambda}}(G_F)$ induces the surjective
map of $l$-adic Lie algebras ${\cal Lie}\, (\rho_{l}(G_F))
\rightarrow {\cal Lie}\, (\rho_{{\Phi}_{\lambda}}(G_F)).
\qed$\enddemo

\proclaim{Lemma 5.22} Let $A/F$ be an abelian variety over $F$ of type I or II
such that $End_{F}\, (A) = End_{{\overline F}}\, (A).$ Then

$$End_{\frak g_{\lambda}}\, (V_{\lambda})\,\cong\,End_{{G_F}}\, (V_{\lambda})\,\,\cong\,\,E_{\lambda} \tag{5.23}$$
$$End_{{\frak g}_{{\Phi}_{\lambda}}}\, (V_{{\Phi}_{{\lambda}}})\,\, \cong
\,\,End_{\Q_l[G_F]}\, (V_{{\Phi}_{{\lambda}}})\,\,\cong\,\,E_{\lambda}. \tag{5.24}$$
\endproclaim
\demo{\sl Proof} By  [F], Theorem 4, the assumption $End_{F}\, (A) =
End_{L}\, (A)$ for any finite extension $L/F,$  Theorem 4.16 (ii), the equality  (5.2) and
Theorem 5.4 we get
$$E_{\lambda} \cong
End_{{E_{\lambda}[G_F]}}\, (V_{\lambda}) \cong End_{{E_{\lambda}[G_L]}}\,
(V_{\lambda}).\tag{5.25}$$ This implies the equality
$$End_{{G_F}}\, (V_{\lambda})= End_{U^{\prime}}\, (V_{\lambda})$$
for any open subgroup $U^{\prime}$ of $G_F.$ Hence, the equality (5.23) follows by Lemma 4.12
(iii). For any $F \subset L \subset {\overline F}$ we have $M_{2,2}(End_{\Q_{l}[G_L]} (V_{l})) =
End_{\Q_{l}[G_L]} (V_{l}^2) = End_{\Q_{l}[G_L]} (V_{l}(A))$ and
$$End_{\Q_{l}[G_L]} (V_{l}(A))\cong\, \prod_{\lambda | l }D_{\lambda}\,\cong \,
\prod_{\lambda | l }M_{2,2}(E_{\lambda}).\tag{5.26}$$ On the other hand
$$\prod_{\lambda | l}E_{\lambda}\,\cong\, \prod_{\lambda | l } End_{E_{\lambda}[G_L]} (V_{\lambda}) \subset
End_{\Q_{l}[G_L]} (V_{l}).\tag{5.27}$$ Hence, comparing the dimensions over $\Q_l$ in (5.26) and
(5.27) we get
$$\prod_{\lambda | l } End_{E_{\lambda}[G_L]}(V_{\lambda})\,
\cong\, End_{\Q_{l}[G_L]} (V_{l}).\tag{5.28}$$
By (5.28) we clearly have
$$\prod_{\lambda | l} End_{\Q_l[G_L]}(V_{{\Phi}_{{\lambda}}}) \subset
End_{\Q_l[G_L]}(V_{l}) \cong \prod_{\lambda | l} E_{\lambda},\tag{5.29}$$
and
$$End_{E_{\lambda}[G_L]}(V_{\lambda}) \subset  End_{\Q_l[G_L]}(V_{{\Phi}_{{\lambda}}}).\tag{5.30}$$
It follows by (5.25), (5.29) and by (5.30) that for any finite field extension $F \subset L$
contained in ${\overline F}$ we have
$$End_{\Q_l[G_L]}(V_{{\Phi}_{{\lambda}}}) \cong
End_{E_{\lambda}[G_L]}(V_{\lambda}) \cong  E_{\lambda}.\tag{5.31}$$
The isomorphisms (5.31) imply
that
$$End_{G_F}(V_{{\Phi}_{{\lambda}}}) \cong End_{U^{\prime}}(V_{{\Phi}_{{\lambda}}})\tag{5.32}$$
for any open subgroup $U^{\prime}$ of $G_F.$ The isomorphism (5.24) follows by (5.32) and Lemma
4.12 (iii).\qed\enddemo
\bigskip

\proclaim{Lemma 5.33} ${\frak g}_{\lambda}^{ss} = sp_{2h} (E_\lambda).$
\endproclaim

\demo{\sl Proof} In the proof we adapt to the current situation the argument from [BGK], Lemma 3.2.
The only thing to check is the miniscule conjecture for the $\lambda$-adic representations
$\rho_{F}\,:\, G_F \rightarrow GL(V_{\lambda}).$ By the work of Pink cf. [P],
Corollary 5.11, we know that ${\frak g}_l^{ss}\otimes \bar\Q_l$ may only have simple factors of
types $A,B,C$ or $D.$ By the semisimplicity of ${\frak g}_l^{ss}$ and Lemma 5.20 the simple
factors of ${\frak g}_{{\Phi}_{\lambda}}^{ss}\otimes \bar\Q_l$ are of the same types. By
Proposition 2.12 and Lemmas 2.21, 2.22, 2.23 we get
$${\frak g}_{{\Phi}_{\lambda}}^{ss} \cong R_{E_{\lambda}/\Q_l}{\frak g}_{\lambda}^{ss}.
\tag{5.34}$$
Since
$${\frak g}_{{\Phi}_{\lambda}}^{ss}\otimes_{\Q_l} {\overline \Q}_l\,\,\cong\,\,
{\frak g}_{\lambda}^{ss} {\otimes}_{E_{\lambda}} E_{\lambda} {\otimes}_{\Q_l} {\overline \Q}
\,\,\cong \,\, \bigoplus_{E_{\lambda}\hookrightarrow {\overline \Q}_l}\,\, {\frak
g}_{\lambda}^{ss}{\otimes}_{E_{\lambda}} {\overline \Q}$$ we see that the simple factors of\,\,
${\frak g}_{\lambda}^{ss}{\otimes}_{E_{\lambda}} {\overline \Q}$\,\, are of types $A,B,C$ or $D.$
The rest of the argument is the same as in the proof of Lemma 3.2 of [BGK].\qed\enddemo

\proclaim{Lemma 5.35} There are natural isomorphisms of  $\Q_l$-algebras.
$$End_{{\frak g}_{{\Phi}_{\lambda}}^{ss}}\, (V_{{\Phi}_{{\lambda}}}) \cong
End_{{\frak g}_{\lambda}^{ss}}\, (V_{\lambda}) \cong  E_{\lambda} \tag{5.36}$$
\endproclaim
\demo{\sl Proof} Since ${\frak g}_{\lambda}$ is reductive and it acts irreducibly on the module
$V_{\lambda}$ (cf.  Lemma 5.33) by [H2], Prop. p. 102 we have:
$${\frak g}_{\lambda} = Z({\frak g}_{\lambda}) \oplus {\frak g}_{\lambda}^{ss} \tag{5.37}$$
and $Z({\frak g}_{\lambda}) =0$ or $Z({\frak g}_{\lambda}) = E_{\lambda}.$
This gives
$$End_{{\frak g}_{\lambda}^{ss}}\, (V_{\lambda}) \,\, =\,\,
End_{{\frak g}_{\lambda}}\, (V_{\lambda}).\tag{5.38}$$ The Weil
restriction functor commutes with the operation of taking the
center of a Lie algebra, hence we get $Z({\frak
g}_{{\Phi}_{{\lambda}}}) = 0$ or $E_{\lambda}$ and by (5.34):
$${\frak g}_{{\Phi}_{{\lambda}}} = Z({\frak g}_{{\Phi}_{{\lambda}}}) \oplus
{\frak g}_{{\Phi}_{{\lambda}}}^{ss}.$$ Since ${\frak
g}_{{\Phi}_{\lambda}} \cong R_{E_{\lambda}/\Q_l}{\frak
g}_{\lambda},$ it is clear that
$$End_{{\frak g}_{{\Phi}_{\lambda}}^{ss}}\, (V_{{\Phi}_{{\lambda}}})\,\, = \,\,
End_{{\frak g}_{{\Phi}_{\lambda}}}\, (V_{{\Phi}_{{\lambda}}}).$$
The lemma follows now from Lemma 5.22. \qed\enddemo

\proclaim{Proposition 5.39} There is an equality of Lie algebras:
$${\frak g}_{l}^{ss}=
\bigoplus_{\lambda |l}\,\, {\frak g}_{{\Phi}_{\lambda}}^{ss} \tag{5.40}$$
\endproclaim

\demo{\sl Proof}
Put ${\overline V_{l}} = V_{l} \otimes_{\Q_l} {\overline \Q}_l,$\quad
${\overline V_{\lambda}} = V_{\lambda} \otimes_{E_{\lambda}} {\overline \Q}_{l},$\quad
${\overline {\frak g}_{l}^{ss}} = {\frak g}_{l}^{ss} \otimes_{\Q_l} {\overline \Q}_{l},$\quad
${\overline {\frak g}_{{\Phi}_{\lambda}}^{ss}} =
{\frak g}_{{\Phi}_{\lambda}}^{ss} \otimes_{\Q_l} {\overline \Q}_{l}.$
By (5.34) we get
$${\overline {\frak g}_{{\Phi}_{\lambda}}^{ss}} \,\,\cong\,\,
{\frak g}_{\lambda}^{ss} \otimes_{E_{\lambda}} E_{\lambda}
\otimes_{\Q_l} {\overline \Q}_l \,\,\cong\,\,
\prod_{E_{\lambda}\hookrightarrow {\overline \Q}_l}
{\frak g}_{\lambda}^{ss} \otimes_{E_{\lambda}} {\overline \Q}_l
\,\, \cong\,\, \prod_{E_{\lambda}\hookrightarrow {\overline \Q}_l}\, sp\, ({\overline V_{\lambda}}) \tag{5.41}$$
By Corollary 1.2.2 of [C1] we have ${\frak g}_{l}=\Q_l \,\oplus\,
{\frak g}_{l}^{ss}$, hence
$$End_{{\frak g}_{l}^{ss}}\, (V_{l}(A)) \,\, = \,\, End_{{\frak g}_{l}}\,
(V_{l}(A)).$$
By Lemmas 5.20 and 5.35
$$\prod_{\lambda | l } E_{\lambda} \,\,\cong\,\,
\prod_{\lambda | l } End_{{\frak g}_{{\Phi}_{\lambda}}^{ss}} (V_{{\Phi}_{\lambda}})\,\, \cong\,\,
\prod_{\lambda | l } End_{{\frak g}_{l}^{ss}} (V_{{\Phi}_{\lambda}})\subset
End_{{\frak g}_{l}^{ss}} (V_{l}).\tag{5.42}$$
But by assumption on $l$ and (5.42)
$$ \prod_{\lambda | l } D_{\lambda} \,\, \cong\,\,
\prod_{\lambda | l } M_{2,2}( E_{\lambda}) \,\, \cong\,\,  M_{2,2}(\prod_{\lambda | l } E_{\lambda})
\subset M_{2,2}(End_{{\frak g}_{l}^{ss}} (V_{l})) \,\,=$$
$$=\,\,  End_{{\frak g}_{l}^{ss}} (V_{l}(A))
\,\, =\,\,  End_{{\frak g}_{l}} (V_{l}(A)) \,\, \cong \,\, \prod_{\lambda | l } D_{\lambda}. \tag{5.43}$$
Comparing dimensions in (5.43) we get
$$End_{{\frak g}_{l}^{ss}} (V_{l}) \,\, \cong\,\, \prod_{\lambda | l } E_{\lambda}. \tag{5.44}$$
Hence we get
$$End_{{\overline {\frak g}_{l}^{ss}}}\, ({\overline V_{l}}) \,\, \cong\,\,
End_{{\frak g}_{l}^{ss}}\, (V_{l}) \otimes_{\Q_l} {\overline \Q}_l \,\, \cong\,\,
\prod_{\lambda | l}
E_{\lambda} \otimes_{\Q_l} \, {\overline \Q}_l \,\, \cong \,\,
\prod_{\lambda | l} \prod_{E_{\lambda}\hookrightarrow {\overline \Q}_l}\,\, {\overline \Q}_l. \tag{5.45}$$
$$End_{{\overline \Q}_l[G_F]}\, ({\overline V_{\lambda}})\,\,\cong\,\,
End_{E_{\lambda}[G_F]}\, (V_{\lambda})\otimes_{E_{\lambda}}{\overline \Q}_l
\,\, \cong\,\,   E_{\lambda} \otimes_{E_{\lambda}}
{\overline \Q}_l \,\, \cong\,\,   {\overline \Q}_l. \tag{5.46}$$
$${\overline V_{l}} \,\,\cong\,\, \bigoplus_{\lambda | l}\,\,
V_{\lambda}\otimes_{\Q_l} \, {\overline \Q}_l \,\,\cong\,\,
\bigoplus_{\lambda | l} \bigoplus_{E_{\lambda}\hookrightarrow {\overline \Q}_{l}}
{\overline V_{\lambda}}. \tag{5.47}$$
By (5.21) the map of Lie algebras
${\overline {\frak g}_{l}^{ss}} \rightarrow {\overline {\frak g}_{{\Phi}_{\lambda}}^{ss}}$
is surjective.
Isomorphisms (5.45), (5.46) and (5.47) show that the simple
${\overline {\frak g}_{l}^{ss}}$ modules
${\frak g}_{\lambda}^{ss} \otimes_{E_{\lambda}} {\overline \Q}_l,$
for all $\lambda | l$ and all $E_{\lambda}\hookrightarrow {\overline \Q}_l,$
are pairwise nonisomorphic submodules of ${\overline {\frak g}_{l}^{ss}}$.
Hence by [H2], Theorem on page 23
$$\bigoplus_{\lambda | l} \bigoplus_{E_{\lambda}\hookrightarrow {\overline \Q}_l}
{\frak g}_{\lambda}^{ss} \otimes_{E_{\lambda}} {\overline \Q}_l
\,\,\subset\,\, {\overline {\frak g}_{l}^{ss}}.\tag{5.48}$$
Tensoring (5.17) with  ${\overline \Q}_l$ and comparing with (5.48) we get
$$\bigoplus_{\lambda | l} \bigoplus_{E_{\lambda}\hookrightarrow {\overline \Q}_l}
{\frak g}_{\lambda}^{ss} \otimes_{E_{\lambda}} {\overline \Q}_l
\,\,\cong \,\, {\overline {\frak g}_{l}}^{ss}.\tag{5.49}$$
Hence for dimensional reasons (5.17), (5.41) and (5.49) imply (5.40).
\qed\enddemo

\proclaim{Corollary 5.50} The representations ${\rho}_{{\Phi}_{\lambda}},$ for $\lambda|l$ are
pairwise nonisomorphic. The representations of the Lie algebra ${\frak g}_{l}^{ss}$ on
$V_{{\Phi}_\lambda}$ are pairwise nonisomorphic over $\Q_l.$
\endproclaim
\demo{\sl Proof} It follows by Lemmas 5.20 and 5.22 and equalities (5.8), (5.36), (5.44).
\qed\enddemo
\bigskip

\proclaim{Corollary 5.51} There is an equality of ranks of group
schemes over $\Q_l$:
$${rank}\,\, (G_{l}^{alg})^{\prime} \,\, =
{rank} \,\, \prod_{\lambda | l} \,\,R_{E_{\lambda}/\Q_l}
(Sp_{2h}/E_{\lambda}).
\tag{5.52}$$
\endproclaim

\demo{\sl Proof} The Corollary follows by Lemma 5.33, equality (5.40), the isomorphism (5.34)
and Lemma 2.21.
\qed\enddemo
\medskip
Taking into account (4.10), (4.11) and Remark 5.13 we get:

$$G(l)^{alg}\,\, \subset \,\,
\prod_{\lambda | l}
R_{k_{\lambda}/\F_l}(GSp_{A_{\lambda}[\lambda]}) \,\, \cong \,\,
\prod_{\lambda | l} R_{k_{\lambda}/\F_l} (GSp_{2h}) \tag{5.53}$$
$$G_{l}^{alg}\,\, \subset \,\,
\prod_{\lambda | l} R_{E_{\lambda}/\Q_l} (GSp_{V_{\lambda}}) \,\,
\cong \,\, \prod_{\lambda | l} R_{E_{\lambda}/\Q_l} (GSp_{2h}).
\tag{5.54}$$

\bigskip

\noindent \subhead  6. Computation of the images of the Galois representations $\rho_l$ and
${\overline\rho}_l$ \endsubhead
\medskip

\noindent In this section we explicitly compute the images of the
l-adic representations induced by the action of the absolute
Galois group on the Tate module of a large class of abelian
varieties of types I and II described in the definition below.

\proclaim{Definition of class $\Cal A$} We say that an abelian
variety $A/F,$ defined over a number field $F,$ is of class ${\cal
A},$ if the following conditions hold:

\roster \item"{(i)}" $A$ is a simple, principally polarized
abelian variety of dimension $g$
 \item"{(ii)}" ${\cal R}=End_{\bar F} (A)=End_{F}(A)$ and the endomorphism
 algebra $D = {\cal R} \otimes_{\Z} \Q,$ is of type I or II in the Albert
list of the  division algebras with involution  cf. [Mu], p. 201
 \item"{(iii)}" the field $F$ is such that for
every $l$ the Zariski closure $G_l^{alg}$ of $\rho_l(G_F)$ in $GL_{2g}/\Q_l$ is a connected
algebraic group
 \item"{(iv)}" $ g =  h e d,$ where $h$ is an odd integer,
 $e = [E: Q]$ is the degree of the center $E$ of $D$ and
 $d^2 = [D:\, E].$
\endroster
\endproclaim

 \noindent Let $L$ be a local
field with the ring of integers ${\cal O}_L$ with maximal ideal ${\frak m}_L = {\frak m}$ and the
residue field $ k = {\cal O}_L/{\frak m}.$

\proclaim{Lemma 6.1}
Let
$$ \diagram {\phantom{\Big|}
{\Cal G}_{1}\phantom{\Big|} }   \ar@{^{(}->}[r]& {\phantom{\Big|}
{\Cal G}_{2}\phantom{\Big|} }
\enddiagram
 \tag{6.2}$$
be an injection of two smooth, reductive group schemes over ${\cal O}_L.$
Let
$$ \diagram {\phantom{\Big|}
G_{1}\phantom{\Big|} }   \ar@{^{(}->}[r]& {\phantom{\Big|}
G_{2}\phantom{\Big|} }
\enddiagram
 \tag{6.3}$$
be the base change to $L$ of the arrow (6.2) and let
$$ \diagram {\phantom{\Big|}
G_{1}({\frak m})\phantom{\Big|} }   \ar@{^{(}->}[r]& {\phantom{\Big|} G_{2}({\frak
m})\phantom{\Big|} }
\enddiagram
 \tag{6.4}$$
be the base change to $k$ of the arrow (6.2). If $rank \, G_1 \, = \, rank\, G_2$\,\, then \,\,
$rank \, G_1({\frak m}) \, = \, rank\, G_2({\frak m}).$
\endproclaim
\demo{\sl Proof}
By [SGA3, Th. 2.5 p. 12] applied to the special point of the scheme
$spec\, {\cal O}_L$ there exists an {\' e}tale neighborhood $S' \rightarrow spec \, {\cal O}_L$
of the geometric point over the special point such that the group schemes
${\cal G}_{1, S'}\, =\, {\cal G}_{1}
\times_{spec\, {\cal O}_L}  S'$  and ${\cal G}_{2, S'}\, =\, {\cal G}_{2}
\times_{spec\, {\cal O}_L}  S'$
have  maximal tori ${\cal T_{1, S'}}$ and ${\cal T_{1, S'}}$ respectively.
By [SGA3] XXII, Th. 6.2.8 p. 260 we observe (we do not need it here but in the Corollary
6.6 below) that
$({\cal G}_{i, S'})' \cap {\cal T_{i, S'}}$ is a maximal torus of
$({\cal G}_{i, S'})'.$ By the definition of a maximal torus
and by
[SGA3] XIX, Th. 2.5, p. 12 applied to the special point of
 $spec\, {\cal O}_L,$ we obtain that the special and generic
fibers of each scheme ${\cal G}_{i, S'}$ have the same rank.
But clearly the generic (resp. special) fibers of
schemes ${\cal G}_{i, S'}$ and ${\cal G}_{i}$ have the same rank for $i = 1, 2.$
Hence going around  the diagram
$$ \diagram {\phantom{\Big|}
G_{1}\phantom{\Big|} }\ar@{_{(}->}[d]
\ar@{^{(}->}[r]& {\phantom{\Big|} G_{2}\phantom{\Big|} }\ar@{_{(}->}[d] \\
{\phantom{\Big|} {\Cal G}_{1}\phantom{\Big|}} \ar@{^{(}->}[r] &
{\phantom{\Big|}{\Cal G}_{2}\phantom{\Big|} }\\
{\phantom{\Big|} G_{1}({\frak m})\phantom{\Big|}} \ar@{^{(}->}[r] \ar@{^{(}->}[u]&
{\phantom{\Big|} G_{2}({\frak m})\phantom{\Big|} }\ar@{^{(}->}[u]\\
\enddiagram
 \tag{6.5}$$

\noindent and taking into account the assumptions that the ranks of the upper corners are the same
we get \, $rank \, G_1({\frak m}) \, = \, rank\, G_2({\frak m}).$ \qed\enddemo

\proclaim{Theorem 6.6} Let $A/F$ be an abelian variety of class
${\cal A}.$ Then for all $l \gg 0,$ we have equalitiy of ranks of
group schemes over $\F_l$:
$${rank}\,\, (G(l)^{alg})^{\prime} \,\, =
{rank} \,\, \prod_{\lambda | l} \,\,R_{k_{\lambda}/\F_l}
(Sp_{2h})
\tag{6.7}$$
\endproclaim
\demo{\sl Proof} By [LP1] Prop.1.3 and by [Wi], Th.1 and 2.1, for $l \gg 0$
the group scheme ${\cal G}_{l}^{alg}$ over $spec\, \Z_l$ is smooth and reductive.
For such an $l$ the structure
morphism $({\cal G}_{l}^{alg})^{\prime} \rightarrow spec\, \Z_l$
is the base
change of the smooth morphism ${\cal G}_{l}^{alg} \rightarrow
D_{\Z_l}(D_{\Z_l}({\cal G}_{l}^{alg}))$ via the unit section
of $D_{\Z_l}(D_{\Z_l}({\cal G}_{l}^{alg})),$ see [SGA3] XXII, Th. 6.2.1, p. 256.
Hence, the group scheme $({\cal G}_{l}^{alg})^{\prime}$ is also smooth over
$\Z_l.$ By [SGA3] {\sl loc. cit}, the group scheme
$({\cal G}_{l}^{alg})^{\prime}$  is semisimple.
We finish the proof by taking $L = \Q_l,$ ${\cal G}_1 =
({\cal G}_{l}^{alg})^{\prime},$ ${\cal G}_2 =
\prod_{\lambda | l} \,\, R_{{\cal O}_{\lambda}/\Z_l} (Sp_{2h})$ in Lemma 6.1 and
applying Corollary 5.51.
\qed\enddemo

\medskip
\remark{\bf Remark 6.8} If $G$ is a group scheme over $S_0$ then the derived subgroup $G^{\prime}$
is defined as the kernel of the natural map
$$G \,\rightarrow D_{S_0}(D_{S_0}(G))$$ [V], [SGA3]. Since this map is consistent with
the base change, we see that for any scheme $S$ over $S_0$ we get
$$G^{\prime}\times_{S_0} S \, = \, (G \times_{S_0} S)^{\prime}.$$
\endremark

\proclaim{Theorem 6.9}  Let $A/F$ be an abelian variety of class
${\cal A}.$ Then for all $l\gg 0,$ we have equalities of group
schemes:

$$(G_{l}^{alg})^{\prime} \,\, = \,\, \prod_{\lambda | l} \,\,R_{E_{\lambda}/\Q_l}
(Sp_{2h})
\tag{6.10}$$
$$(G(l)^{alg})^{\prime} \,\, = \,\, \prod_{\lambda | l} \,\,R_{k_{\lambda}/\F_l}
(Sp_{2h})
\tag{6.11}$$
\endproclaim
\demo{\sl Proof} The proof is similar to the proof of Lemma 3.4 of
[BGK]. We prove the equality (6.11). The proof of the equality
(6.10) is analogous. Let $$\underline{\rho}_{l}\,:\, G(l)^{alg}
\rightarrow  \, GL_{2g}$$ denote the representation induced by the
inclusion $G(l)^{alg}\,\subset\, GL_{2g}.$ By the result of Faltings
cf. [Fa], the representation $\underline{\rho}_{l}$ is semisimple
and the commutant of $\underline{\rho}_{l}(G(l)^{alg})$ in the
matrix ring $M_{2g,2g}$ is $End_{\bar F} (A) \otimes_{\Z} \F_l.$ The
representation $\underline{\rho}_{l}$ factors through the imbedding
(5.53). Projecting onto the $\lambda$ component in (5.53) we obtain
the representation
$$\underline{\rho}_{\Phi_{\lambda}}\,:\,
G(l)^{alg} \rightarrow  \, R_{k_{\lambda}/\F_l} ( GSp_{A[\lambda]} ) \,\cong\,
R_{k_{\lambda}/\F_l} ( GSp_{2h} ). \tag{6.12}$$
This map corresponds to the map
$$ G(l)^{alg} \otimes_{\F_l} k_{\lambda}  \rightarrow  \, GSp_{2h}. \tag{6.13}$$
By Remark 6.8 restriction of the the map (6.13) to the derived subgroups
gives the following map:
$$ (G(l)^{alg})^{\prime} \otimes_{\F_l} k_{\lambda}  \rightarrow  \, Sp_{2h} \tag{6.14}$$
which in turn gives the representation
$$\underline{\rho}_{\Phi_{\lambda}}\,:\,
(G(l)^{alg})^{\prime} \rightarrow  \, R_{k_{\lambda}/\F_l} (
Sp_{2h} ).$$ Now by (5.3) we have the natural isomorphisms:
$$\prod_{k_{\lambda}\hookrightarrow {\overline \F}_l} {\overline \F}_l \,\,\cong\,\,
k_{\lambda}\otimes_{\F_{l}} {\overline \F}_l \,\,\cong\,\,
End_{k_{\lambda}\otimes_{\F_{l}} {\overline \F}_l [G_F]} (A_{\lambda}[\lambda]
\otimes_{\F_{l}} {\overline \F}_l) \,\,\cong $$
$$ \cong\,\,
End_{k_{\lambda}\otimes_{\F_{l}} {\overline \F}_l [G_F]} (A_{\lambda}[\lambda]
\otimes_{k_{\lambda}} k_{\lambda}
\otimes_{\F_{l}} {\overline \F}_l)\,\, \cong $$
$$\cong \,\,
\prod_{k_{\lambda}\hookrightarrow {\overline \F}_l}
End_{{\overline \F}_l [G_F]} (A_{\lambda}[\lambda]
\otimes_{k_{\lambda}} {\overline \F}_l).\tag{6.15}$$ 
Note that $Z(Sp_{2h}) \cong \mu_{2}$ and this isomorphism holds over any field of
definition. The isomorphisms (6.15) imply by the Schur's Lemma:
$$\underline{\rho}_{\Phi_{\lambda}} (Z((G(l)^{alg})^{\prime}))
\subset R_{k_{\lambda}/\F_l}( \mu_{2}).$$ Hence
$$Z((G(l)^{alg})^{\prime}) \subset \prod_{\lambda | l } R_{k_{\lambda}/\F_l}( \mu_{2})
= Z(\prod_{\lambda | l } R_{k_{\lambda}/\F_l}( Sp_{2h})).$$
Observe that both groups $(G(l)^{alg})^{\prime}$ and
$\prod_{\lambda | l } R_{k_{\lambda}/\F_l}( Sp_{2h})$ are
reductive. Now the proof is finished in the same way as the proof
of Lemma 3.4 in [BGK]. \qed\enddemo

\noindent \proclaim{Theorem 6.16}  Let $A/F$ be an abelian variety of class ${\cal A}.$ Then for
$l \gg 0,$ we have:

\noindent
$${\overline{\rho_l}}(G_F^{\prime})
\,\, = \,\, \prod_{\lambda | l} \,\, Sp_{2h}(k_{\lambda})
\,\, =\,\, Sp_{2h}({\cal O}_{E}/ l {\cal O}_{E}),
\tag{6.17}$$

$$\rho_l (\overline{G_F^{\prime}})
\,\, = \,\, \prod_{\lambda | l} \,\, Sp_{2h}({\cal O}_{\lambda}) \,\, =\,\, Sp_{2h}({\cal O}_{E}
\otimes_{\Z} \Z_l),\tag{6.18}$$ where ${\overline{\rho_l}}$ is the representation ${\rho_l}$
$mod\,\,\, l$ and $\overline{G_F^{\prime}}$ is the closure of the commutator subgroup
$G_F^{\prime}\subset G_F$ computed with respect to the natural profinite topology of $G_F.$
\endproclaim

\demo{\sl Proof} To prove the equality (6.17), note that the group
scheme $\prod_{\lambda | l } R_{k_{\lambda}/\F_l}( Sp_{2h})$ is
simply connected, since its base change to ${\overline \F}_l$ is
$\prod_{\lambda | l } \prod_{k_{\lambda}\hookrightarrow {\overline
\F}_l} Sp_{2h}/{\overline \F}_l,$ which is clearly simply
connected. From now on the argument is the same as in the proof of
Theorem 3.5 in [BGK]. Namely: it follows by (6.11) that
$(G(l)^{alg})^{\prime}$ is simply connected. So
$(G(l)^{alg})^{\prime}(\F_l) = (G(l)^{alg})^{\prime}(\F_l)_u.$
Hence, by a theorem of Serre (cf. [Wi], Th.4) we get
$$(G(l)^{alg})^{\prime}(\F_l) \,\,\subset\,\,
({\overline{\rho_l}}(G_F))^{\prime} \,\, = \,\,
{\overline{\rho_l}}(G_{F}^{\prime}).$$ On the other hand, by
definition of the group $G(l)^{alg},$ it is clear that
$${\overline{\rho_l}}(G_F^{\prime}) = ({\overline{\rho_l}}(G_F))^{\prime}\,\, \subset\,\,
(G(l)^{alg})^{\prime}(\F_l).$$ As for the second equality in  (6.18)
we have
$$\rho_l\bigl(\,\overline{G_F^{\prime}}\,\bigr) \,\, = \,\,
{\overline{({\rho_l}(G_F))^{\prime}}} \,\, \subset \,\,
\prod_{\lambda | l} \,\, Sp_{2h}({\cal O}_{\lambda}), \tag{6.19}$$
where ${\overline{({\rho_l}(G_F))^{\prime}}}$ denotes the closure
of $({\rho_l}(G_F))^{\prime}$  in the natural ($\lambda$-adic in
each factor) topology of the group $\prod_{\lambda | l} \,\,
Sp_{2h}({\cal O}_{\lambda}).$  Using equality (6.17) and Lemma 6.20
stated below, applied to $X =
{\overline{({\rho_l}(G_F))^{\prime}}},$ we finish the proof. \qed
\enddemo

\proclaim{Lemma 6.20} Let $X$ be a closed subgroup in
$\prod_{\lambda | l} Sp_{2h}({\cal O}_{\lambda})$ such that its
image via the reduction map
$$\prod_{\lambda | l} Sp_{2h}({\cal O}_{\lambda}) \rightarrow
\prod_{\lambda | l} Sp_{2h}(k_{\lambda})$$
is all of $\prod_{\lambda | l} Sp_{2h}(k_{\lambda}).$ Then
$X = \prod_{\lambda | l} Sp_{2h}({\cal O}_{\lambda}).$
\endproclaim
\demo{\sl Proof} The proof is similar to the proof of Lemma 3 in [Se]
chapter IV, 3.4.
\qed
\enddemo
\bigskip

\noindent

\subhead  7. Applications to classical conjectures
\endsubhead
\bigskip

\noindent
Choose an imbedding of $F$ into the field of complex numbers ${\Bbb C}.$ Let $V =
H^1(A({\Bbb C}), \Q)$ be the singular cohomology group with rational coefficients. Consider the
Hodge decomposition
$$V\otimes_{\Q} {\Bbb C}=H^{1,0}\oplus H^{0,1},$$
where \,\,$H^{p,q} = H^p(A;\, \Omega_{A/\C}^q)$ and ${\overline{H^{p,q}}}= H^{q,p}.$ Observe that
$H^{p,q}$ are invariant subspaces with respect to $D = End_{{\overline F}} (A) \otimes \Q$ action
on $V \otimes_{\Q} \C.$ Hence, in particular $H^{p,q}$ are $E$-vector spaces. Let
$$\psi\, :\, V \times V\, \rightarrow \Q$$
be the $\Q$-bilinear,  nondegenerate, alternating form coming from the Riemann form of $A.$ Since
$A$ has a principal polarization by assumption, the form $\psi$ is given by the standard matrix
$$J =
\pmatrix  0 & I_g\\
                - I_g& 0\\
\endpmatrix.
$$
Define the cocharacter
$$\mu_{\infty}:{\Bbb G}_{m}({\Bbb C})\rightarrow GL(V\otimes_{\Q}
{\Bbb C})=GL_{2g}({\Bbb C})$$ such that, for any $z\in {\Bbb
C}^{\times},$ the automorphism $\mu_{\infty}(z)$ is the multiplication by
$z$ on $H^{1,0}$ and the identity on $H^{0,1}.$

\proclaim{Definition 7.1} The Mumford-Tate group of the abelian variety $A/F$ is the smallest
algebraic subgroup $MT(A) \subset GL_{2g},$ defined over $\Q,$ such that $MT(A)({\Bbb C})$
contains the image of $\mu_{\infty}.$ The Hodge group $H(A)$ is by definition the connected
component of the identity in $MT(A) \cap SL_{V} \cong MT(A) \cap SL_{2g}.$
\endproclaim

\noindent We refer the reader to [D] for an excellent exposition
on the Mumford-Tate group. In particular, $MT(A)$ is a reductive
group {\sl loc. cit.} Since, by definitions
$$\mu_{\infty}(\C^{\times}) \, \subset GSp_{(V,\, \psi)} (\C) \cong GSp_{2g}(\C),$$
it follows that the group $MT(A)$ is a reductive subgroup of the
group of symplectic similitudes $GSp_{(V,\, \psi)} \cong GSp_{2g}$
and that
$$H(A)\, \subset Sp_{(V,\, \psi)} \cong Sp_{2g}.\tag{7.2}$$

\remark{\bf Remark 7.3} Let $V$ be a finite dimensional vector space
over a field $K$ such that it is also an $R$-module for a
$K$-algebra $R.$ Let $G$ be a $K$-group subscheme of $GL_V.$ Then by
the symbol $C_R(G)$ we will denote the commutant of $R$ in $G.$ The
symbol $C_R^0(G)$ will denote the connected component of identity in
$C_R(G).$ Let $\beta \,:\, V \times V \rightarrow K$ be a bilinear
form and let $G_{(V, \beta)} \subset GL_V$ be the subscheme of
$GL_V$ of all  isometries with respect to the bilinear form $\beta.$
It is easy to check that \,\, $C_{R}(G_{(V, \beta)}) \otimes_K L
\cong C_{R \otimes_K L}(G_{(V \otimes_K L,\,\beta \otimes_K L)}).$
Note that $MT(A) \subset C_{D}(GSp_{(V,\, \psi)})$ by  definitions.
\endremark

\noindent \proclaim{Definition 7.4} The algebraic group $L(A) = C_D^0(Sp_{(V,\, \psi)})$ is called
the {\it Lefschetz group} of a principally polarized abelian variety $A.$ Note that the group
$L(A)$ does not depend on the form $\psi$ cf. [R2].
\endproclaim

\noindent By [D], Sublemma 4.7, there is a unique $E$-bilinear, nondegenerate, alternating pairing
$$\phi\, :\, V \times V\, \rightarrow E$$
such that $Tr_{E/\Q} (\phi) = \psi.$ Taking into account that the actions of $H(A)$ and $L(A)$ on
$V$ commute with the $E$-structure, we get
$$H(A) \, \subset \, L(A) \,\subset \,R_{E/\Q} Sp_{(V,\,
\phi)} \, \subset \, Sp_{(V,\, \psi)}.\tag{7.5}$$ But $R_{E/\Q}
(Sp_{(V,\, \phi)}) = C_E(Sp_{(V,\, \psi)})$ hence $C_D(R_{E/\Q}
(Sp_{(V,\, \phi)})) = C_D(Sp_{(V,\, \psi)})$ so
$$H(A) \subset L(A) = C_D^0(R_{E/\Q} (Sp_{(V,\, \phi)})) \subset
C_D(R_{E/\Q} (Sp_{(V,\, \phi)})).\tag{7.6}$$
\medskip

\noindent \proclaim{Definition 7.7} If $L/\Q$ is a field extension
of $\Q$ we put
$$MT(A)_{L} := MT(A) \otimes_{\Q} L,\quad
H(A)_{L} := H(A) \otimes_{\Q} L,\quad L(A)_{L} := L(A) \otimes_{\Q} L.$$
\endproclaim
\medskip

\noindent \proclaim{Conjecture 7.8 (Mumford-Tate cf. [Se5],
C.3.1)} If $A/F$ is an abelian variety over a number field $F$,
then for any prime number $l$
$$(G_l^{alg})^o= MT(A)_{\Q_l},\tag{7.9}$$
where $(G_l^{alg})^o$ denotes the connected component of the
identity.
\endproclaim

\proclaim{Theorem 7.10 (Deligne [D], I, Prop. 6.2)} If $A/F$ is an
abelian variety over a number field $F$ and $l$ is a prime number,
then
$$(G_l^{alg})^o\subset MT(A)_{\Q_l}.\tag{7.11}$$
\endproclaim

\proclaim{Theorem 7.12} The Mumford-Tate conjecture holds true for
abelian varieties of class ${\cal A}$ defined in the beginning of
Section 6.
\endproclaim

\demo{Proof} By [LP1], Theorem 4.3, it is enough to verify (7.9)
for a single prime $l$ only. We use the equality (6.10) for a big
enough prime $l.$ The proof goes similarly to the proof of Theorem
3.6 in [BGK]. In the proof we will make some additional
computations, which provide an extra information on the Hodge
group $H(A).$ The Hodge group $H(A)$ is semisimple (cf. [G], Prop.
B.63) and the center of $MT(A)$ is $\G_{m}$ (cf. [G], Cor. B.59).
Since $MT(A) =\G_{m} H(A),$ we get
$$(MT(A)_{\Q_l})^{\prime} = (H(A)_{\Q_l})^{\prime} = H(A)_{\Q_{l}}. \tag{7.13}$$
By (7.11), (7.13)  and (6.10)
$$\prod_{\lambda | l} R_{E_{\lambda}/\Q_l} ( Sp_{(W_{\lambda},
\psi_{\lambda}^{0})}) \cong \prod_{\lambda | l}
R_{E_{\lambda}/\Q_l} ( Sp_{2h}) \subset H(A)_{\Q_l}. \tag{7.14}$$
On the other hand by (7.6)
$$H(A)_{\Q_l}
\subset L(A)_{\Q_l}\subset C_D(R_{E/\Q} (Sp_{(V,\,
\phi)}))\otimes_{\Q} \Q_l . \tag{7.15}$$ Since $R_{E/\Q}
(Sp_{(V,\, \phi)}) = C_E(Sp_{(V,\, \psi)}),$ by Remark 7.3,  formulae
(7.14) and  (7.15) we get:
$$\prod_{\lambda | l} R_{E_{\lambda}/\Q_l} ( Sp_{(W_{\lambda}, \psi_{\lambda}^{0})})
\subset \prod_{\lambda | l} C_{D_{\lambda}}(R_{E_{\lambda}/\Q_l}
(Sp_{(V_{\lambda}(A),\, \psi_{\lambda}^{0})})). \tag{7.16}$$ For $A$
of type I, $D_{\lambda} = E_{\lambda}$ and $V_{\lambda}(A) =
W_{\lambda}$ hence, trivially, the inclusion (7.16) is an equality.
Assume that $A$ is of type II. Since $V_{\lambda}(A) = W_{\lambda}
\oplus W_{\lambda}$ and $D_{\lambda} = M_{2,2} (E_{\lambda}),$
evaluating both sides of the inclusion (7.16) on the ${\overline
\Q}_l$-points, we get equality with both sides equal to
$$\prod_{\lambda | l} \prod_{E_{\lambda} \hookrightarrow {\overline \Q}_l} ( Sp_{(W_{\lambda},
\phi_{\lambda} |_{W_{\lambda}})})({\overline \Q}_l )$$ which is an
irreducible algebraic variety over ${\overline \Q}_l.$ Then we use
Prop. II, 2.6 and Prop. II, 4.10 of [H] in order to conclude that
the  groups $H(A)_{{\overline \Q}_l},$ $L(A)_{{\overline \Q}_l}$
and $C_D(R_{E/\Q} (Sp_{(V,\, \phi)}))\otimes_{\Q} {\overline
\Q}_l$ are connected. Hence all the groups $H(A),$ $L(A)$ and
$C_D(R_{E/\Q} (Sp_{(V,\, \phi)}))$ are  connected, and we have
$$\prod_{\lambda | l} R_{E_{\lambda}/\Q_l} ( Sp_{(W_{\lambda},
\phi_{\lambda} |_{W_{\lambda}})}) \cong \prod_{\lambda | l} R_{E_{\lambda}/\Q_l} ( Sp_{2h}) =
\tag{7.17}$$
$$= H(A)_{{\Q}_l} =  L(A)_{{\Q}_l} = C_D(R_{E/\Q} (Sp_{(V,\, \phi)}))\otimes_{\Q} {\Q}_l.$$
By (6.10), (7.17) and [Bo], Corollary 1. p. 702 we get
$$MT(A)_{\Q_l} = \G_{m} H(A)_{\Q_l} = \G_{m} (G_{l}^{alg})^{\prime} \subset G_{l}^{alg}. \tag{7.18}$$
The Theorem follows by (7.11) and (7.18). \qed\enddemo

\proclaim{Corollary 7.19} If $A$ is an abelian variety of class
${\cal A},$ then
$$H(A)_{\Q} = L(A)_{\Q} = C_D(R_{E/\Q} (Sp_{(V,\, \phi)})) =
C_D(Sp_{(V,\, \psi)}) .\tag{7.20}$$
\endproclaim
\demo{Proof} Taking Lie algebras of groups in (7.17) we deduce by a simple dimension argument that
$${\cal Lie}\,  H(A) =  {\cal Lie} \, L(A) = {\cal Lie} \, C_D(R_{E/\Q} (Sp_{(V,\, \phi)})).
\tag{7.21}$$ In the proof of Theorem 7.12 we have showed that the
groups $H(A),$ $L(A)$ and $C_D(R_{E/\Q} (Sp_{(V,\, \phi)}))$ are
connected. Hence, by Theorem p. 87 of [H1] we conclude that
$$H(A) =  L(A) = C_D(R_{E/\Q} (Sp_{(V,\, \phi)})).\tag{7.22} $$\enddemo

\proclaim{Corollary 7.23}  If $A$ is an abelian variety of class ${\cal A},$ then for all $l$:
$$H(A)_{\Q_l} = \prod_{\lambda | l}
C_{D_{\lambda}}(R_{E_{\lambda}/\Q_l} (Sp_{(V_{\lambda}(A),\, \phi
\otimes_{\Q} E_{\lambda})})). \tag{7.24}$$
In particular, for $l
\gg 0$ we get
$$H(A)_{\Q_l} = \prod_{\lambda | l}
R_{E_{\lambda}/\Q_l} (Sp_{(W_{\lambda},\, \phi \otimes_{\Q} E_{\lambda})}).\tag{7.25}$$
\endproclaim
\demo{Proof} Equality (7.24) follows immediately from Corollary 7.19. Equality (7.25) follows then
from (7.17). \qed\enddemo
\bigskip

\noindent We have:
$$H^1(A(\C);\, \R)\, \cong \, V \otimes_{\Q} \R \, \cong\, \bigoplus_{\sigma : E \hookrightarrow \R}
V \otimes_{E, \sigma} \R.$$ Put $V_{\sigma}(A) =  V \otimes_{E,
\sigma} \R$ and let $\phi_{\sigma}$  be the form
$$\phi \otimes_{E, \sigma} \R\, :\,  V_{\sigma}(A) \otimes_{\R} V_{\sigma}(A)\,
\rightarrow \, \R.$$

\proclaim{Lemma 7.26} If $A$ is simple, principally polarized
abelian variety of type II, then for each $\sigma : E
\hookrightarrow \R $ there is an $\R$-vector space $W_{\sigma}(A)$
of dimension ${{g}\over{e}} = {4\,  dim A \over [D;\, \Q]}$ such
that: \roster \item"{(i)}"$V_{\sigma}(A) \cong W_{\sigma}(A)
\oplus W_{\sigma}(A),$ \item"{(ii)}"the restriction of \,\, $\phi
\otimes_{\Q} \R$\,\, to $W_{\sigma}(A)$ gives a nondegenerate,
alternating pairing
$$\psi_{\sigma}\,:\, W_{\sigma}(A) \times W_{\sigma}(A)
\rightarrow \R.$$
\endroster
\endproclaim

\demo{Proof} Using the assumption that $D \otimes_{\Q} \R \cong M_{2,2}(\R)$ the proof is similar
to the proof of Theorem 5.4. \qed\enddemo

\noindent
We put $$ W_{\infty, \sigma}  = \cases V_{\sigma}(A)   & \text{if $A$ is of type I }\\
\\
W_{\sigma}(A)\, ,  & \text{if $A$ is of type II }
\endcases$$
and
$$ \psi_{\sigma}  = \cases \phi_{\sigma}   & \text{if $A$ is of type I }\\
\\
\phi_{\sigma} |_{W_{\sigma}(A)} \, ,  & \text{if $A$ is of type II. }
\endcases$$

\noindent Observe that
$$ dim_{\R}\, W_{\infty, \sigma}  =
\cases {2 g \over e} = {2\,  dim A \over [D;\, \Q]}  & \text{if $A$ is of type I} \\
\\
{g \over e} = {4\,  dim A \over [D;\, \Q]} \, ,  & \text{if $A$ is of type II. }
\endcases$$

\proclaim{Corollary 7.27}  If $A$ is an abelian variety of class
${\cal A},$ then
$$H(A)_{\R} =  L(A)_{\R} = \prod_{\sigma : E \hookrightarrow \R}
Sp_{(W_{\infty, \sigma},\, \psi_{\sigma})}\tag{7.28}$$
$$H(A)_{\C} =  L(A)_{\C} = \prod_{\sigma E \hookrightarrow \R}
Sp_{(W_{\infty, \sigma} \otimes_{\C} \C,\, \psi_{\sigma}  \otimes_{\R} \C)}.\tag{7.29}$$
\endproclaim
\demo{Proof} It follows from Lemma 7.26 and Corollary 7.19. \qed\enddemo
\medskip

\noindent We recall the conjectures of Tate and Hodge in the case of
abelian varieties. See [G], [K] and [T1] for more details.

\proclaim{Conjecture 7.30 (Hodge)} If $A/F$ is a simple abelian
variety over a number field $F,$ then for every $0 \leq p \leq g$
the natural cycle map induces an isomorphism
$$A^p(A) \cong H^{2p}(A(\C);\, \Q) \cap H^{p,p},\tag{7.31}$$
where \,\,$A^p(A)$ is the $\Q$-vector space of codimension $p$
algebraic cycles on $A$ modulo the homological equivalence.
\endproclaim

\proclaim{Conjecture 7.32 (Tate)} If $A/F$ is a simple abelian
variety over a number field $F$ and $l$ is a prime number, then
for every $0 \leq p \leq g$ the cycle map induces an isomorphism:
$$A^p(A) \otimes_{\Q} \Q_l \cong H^{2p}_{et}({\overline A};\, \Q_l(p))^{G_F}\tag{7.33}$$
where  ${\overline A} = A \otimes_{F} {\overline F}.$
\endproclaim
\medskip

\noindent \proclaim{Theorem 7.34} The Hodge conjecture holds true
for abelian varieties of class ${\cal A}.$
\endproclaim

\demo{Proof} By [Mu], Theorem 3.1 the Hodge conjecture follows from the equality (7.20) of
Corollary 7.19. \qed\enddemo
\medskip

\noindent \proclaim{Theorem 7.35} The Tate conjecture holds true
for abelian varieties of class ${\cal A}.$
\endproclaim

\demo{Proof} It is known (see Proposition 8.7 of [C1]) that
Mumford-Tate conjecture implies the equivalence of Tate and Hodge
conjectures. Hence the Tate conjecture follows by Theorems 7.12
and 7.34. \qed\enddemo
\medskip

\noindent \proclaim{Conjecture 7.36 (Lang)} Let $A$ be an abelian variety over a number field $F.$
Then for $l \gg 0$ the group $\rho_{l}(G_F)$ contains the group of all  homotheties
 in $GL_{T_l(A)}(\Z_l).$
\endproclaim

\noindent \proclaim{Theorem 7.37 (Wintenberger [Wi], Cor. 1, p.
5)} Let $A$ be an abelian variety over a number field $F$. The
Lang conjecture holds for such abelian varieties $A$ for which the
Mumford-Tate conjecture holds or if $dim\, A < 5.$
\endproclaim

\noindent \proclaim{Theorem 7.38} The Lang's conjecture holds true
for abelian varieties of class ${\cal A}.$
\endproclaim

\demo{Proof} It follows by Theorem 7.12 and Theorem 7.37. \qed\enddemo
\medskip

\noindent We are going to use Theorem 7.12 and Corollary 7.19 to prove an analogue of the open
image Theorem of Serre cf. [Se8]. We start with the following remark which is a plain
generalization of remark 7.3.

\remark{\bf Remark 7.39} Let $B_1 \subset B_2$ be two commutative
rings with identity. Let $\Lambda$ be a free, finitely generated $B_1$-module 
such that it is also an $R$-module for a $B_1$-algebra $R.$
Let $G$ be a $B_1$-group subscheme of $GL_{\Lambda}.$ Then
$C_R(G)$ will denote the commutant of $R$ in $G.$ The symbol
$C_R^0(G)$ will denote the connected component of identity in
$C_R(G).$ Let $\beta \,:\, \Lambda \times \Lambda \rightarrow B_1$
be a bilinear form and let $G_{(\Lambda, \beta)} \subset
GL_{\Lambda}$ be the subscheme of $GL_{\Lambda}$ of the isometries
with respect to the  form $\beta.$ Then we check that \,\, $C_{R}(G_{({\Lambda},
\beta)}) \otimes_{B_1} B_2 \cong C_{R \otimes_{B_1}
B_2}(G_{(\Lambda \otimes_{B_1} B_2,\,\beta \otimes_{B_1} B_2)}).$
\endremark
\medskip

\noindent
Consider the bilinear form:
$$\psi\,\, :\,\, \Lambda \times \Lambda \, \rightarrow \, \Z\tag{7.40}$$
associated with the variety $A.$ Abusing notation sligthly, we will denote by $\psi$ the Riemann
form $\psi \otimes_{\Z} \Q,$ i.e., we put:
$$\psi\,\, :\,\, V \times V \, \rightarrow \, \Q.$$
Consider the group scheme $C_{{\cal R}}(Sp_{(\Lambda,\, \psi)})$
over $Spec\, \Z.$ Since $C_{{\cal R}}(Sp_{(\Lambda,\, \psi)})
\otimes_{\Z} \Q = C_{D}(Sp_{(V,\, \psi)})$ (see Remark  7.39),
there is an open imbedding in the $l$-adic topology:
$$C_{{\cal R}}(Sp_{(\Lambda,\, \psi)})(\Z_{l}) \subset C_{D}(Sp_{(V,\,
\psi)})(\Q_{l}).\tag{7.41}$$ Note that the form $\psi_l$ of (4.1)
is obtained by tensoring  (7.40) with $\Z_l.$
\medskip

\noindent \proclaim{Theorem 7.42} If $A$ is an abelian variety of class ${\cal A},$ then for every
prime number $l,$ $\rho_{l}(G_F)$ is open in the group
$$C_{{\cal
R}}(GSp_{(\Lambda,\, \psi)})(\Z_l) = C_{{\cal R}\otimes_{\Z} \Z_l}(GSp_{(T_{l}(A),\,
\psi_l)})(\Z_{l}).$$ In addition, for $l \gg 0$ we have:
$$\rho_{l}(\overline{G_F^{\prime}}) \,\, =\,\,
C_{{\cal R}}(Sp_{(\Lambda,\, \psi)})(\Z_l). \tag{7.43}$$
\endproclaim
\demo{Proof} For any ring with identity $R$ the group
$GSp_{2g}(R)$ is generated by subgroups $Sp_{2g}(R)$ and
$$\{
\pmatrix  a I_g & 0\\
                0 & I_g\\
\endpmatrix;\, a \in R^{\times}
\}. $$ One checks easily that the group $\Z_{l}^{\times} Sp_{2g}(\Z_l)$ has index $2$ (index $4$ resp.)
in $GSp_{2g}(\Z_l),$ for $l > 2$ (for $l = 2$ resp.). Here the symbol $\Z_{l}^{\times}$ stands for the
subgroup of homotheties in $GL_{2g}(\Z_l).$ Since  by assumption  $A$ has a principal polarization,
$Sp_{2g}(\Z_l) \cong Sp_{(\Lambda,\, \psi)})(\Z_l).$ By [Bo], Cor. 1. on p. 702,  there is an open
subgroup $U \subset \Z_{l}^{\times}$ such that $U \subset \rho_l(G_F).$ Hence $U \, C_{{\cal
R}}(Sp_{(\Lambda,\, \psi)})(\Z_l) = C_{{\cal R}} (U \, Sp_{(\Lambda,\, \psi)}(\Z_l))$ is an open
subgroup of $C_{{\cal R}}(GSp_{(\Lambda,\, \psi)})(\Z_l) = C_{{\cal R}}(GSp_{(\Lambda,\,
\psi)}(\Z_l)).$ By [Bo], Th. 1, p. 701,  the group $\rho_{l}(G_F)$ is open in $G_{l}^{alg}
(\Q_l).$ By Theorem 7.12, Corollary 7.19 and Remark 7.3
$$U \, C_{{\cal R}}(Sp_{(\Lambda,\, \psi)})(\Z_l) \,\, \subset \,\,
\Q_{l}^{\times}\, C_{D}(Sp_{(V,\, \psi)})(\Q_l)\,\, = $$
$$ = \,\,\G_{m}(\Q_l) H(A)(\Q_l) \subset
MT(A) (\Q_l) = G_{l}^{alg} (\Q_l).\tag{7.44}$$ Hence, $U \, C_{{\cal R}}(Sp_{(\Lambda,\,
\psi)})(\Z_l) \cap \rho_{l}(G_F)$ is open in $U \, C_{{\cal R}}(Sp_{(\Lambda,\, \psi)})(\Z_l)$ and we
get that $\rho_{l}(G_F)$ is open in $C_{{\cal R}}(GSp_{(\Lambda,\, \psi)})(\Z_l).$ Using Remark
7.39 and the universality of the fiber product, we observe that
$$C_{{\cal R}}(Sp_{(\Lambda,\, \psi)})(\Z_{l}) =
C_{{\cal R}\otimes_{\Z} \Z_l}(Sp_{(T_{l}(A),\, \psi_l)})(\Z_{l}).
\tag{7.45}$$ For $l \gg 0$ we get
$$C_{{\cal R}\otimes_{\Z} \Z_l}(Sp_{(T_{l}(A),\, \psi_l)}) \,\,\cong\,\,
C_{{\cal R}\otimes_{\Z} \Z_l}(C_{{\cal O}_{E} \otimes_{\Z} \Z_l} (Sp_{(T_{l}(A),\, \psi_l)})) \,\,
\cong\,\, $$
$$ \cong\,\,  C_{{\cal R}\otimes_{\Z} \Z_l}(\prod_{\lambda | l} R_{{\cal O}_{\lambda}/\Z_l}
(Sp_{(T_{\lambda}(A),\, \psi_{\lambda})})). \tag{7.46}$$ Evaluating the group schemes in (7.46) on
$Spec\, \Z_l$ we get
$$C_{{\cal R}\otimes_{\Z} \Z_l}(Sp_{(T_{l}(A),\, \psi_l)})(\Z_l) \,\,\cong\,\,
C_{{\cal R}\otimes_{\Z} \Z_l}(\prod_{\lambda | l} R_{{\cal O}_{\lambda}/\Z_l}
(Sp_{(T_{\lambda}(A),\, \psi_{\lambda})}))(\Z_l)\,\,\cong $$
$$\cong \,\, \prod_{\lambda | l} C_{{\cal R}_{\lambda}}
Sp_{(T_{\lambda}(A),\, \psi_{\lambda})}({\cal O}_{\lambda})\,\,
\cong \,\, \prod_{\lambda | l} Sp_{(T_{\lambda},\,
\psi_{\lambda})}({\cal O}_{\lambda}) \cong \,\, \prod_{\lambda |
l} Sp_{2h}({\cal O}_{\lambda}). \tag{7.47}$$ Hence by (7.45),
(7.46), (7.47), (6.18) and Theorem 7.38,  we conclude that for
$l\gg 0$ the equality (7.43) holds. \qed\enddemo
\medskip

\noindent \proclaim{Theorem 7.48} If $A$ is an abelian variety of class ${\cal A},$ then for every
prime number $l,$ the group $\rho_{l}(G_F)$ is open in the group ${\cal G}_{l}^{alg}(\Z_l)$ in the
$l$-adic topology.
\endproclaim
\demo{Proof} By Theorem 7.42 the group $\rho_{l}(G_F)$ is open
 in $C_{{\cal R}\otimes_{\Z} \Z_l}(GSp_{(T_{l}(A),\, \psi_l)})(\Z_{l})$
is an open imbedding in the $l$-adic topology, so $\rho_{l}(G_F)$
has a finite index in the group \newline $C_{{\cal R}\otimes_{\Z}
\Z_l}(GSp_{(T_{l}(A),\, \psi_l)})(\Z_{l}).$ By the definition of
${\cal G}_{l}^{alg},$ we have:
 $$\rho_{l}(G_F) \subset {\cal G}_{l}^{alg}(\Z_l) \subset C_{{\cal R}\otimes_{\Z} \Z_l}(GSp_{(T_{l}(A),\,
 \psi_l)})(\Z_{l}).$$
Hence, $\rho_{l}(G_F)$ has a finite index in ${\cal G}_{l}^{alg}(\Z_l),$ and the claim
follows since
\newline
 $ C_{{\cal R}\otimes_{\Z} \Z_l}(GSp_{(T_{l}(A),\,
 \psi_l)})(\Z_{l})$ is a profinite group. \qed\enddemo
\medskip

\noindent {\it Acknowledgements}.\,\,\,  The first author would like to thank G. Faltings for
conversations concerning l-adic representations associated with abelian varieties. The first and
the second authors would like to thank Mathematics Department of Ohio State University and Max
Planck Institut in Bonn for financial support during visits in 2003 and 2004. The second author
thanks the CRM in Barcelona for the hospitality during the visit in June and July 2004. The
research has been partially sponsored by a KBN grant. 

\enddocument